\documentclass[12pt,a4paper]{article}

\usepackage[table]{xcolor}
\usepackage{slashed}
\usepackage{xkeyval}
\usepackage{url,amsmath,amssymb,latexsym,pstricks,mathrsfs,comment,amsthm,graphicx,tikz,tikz-cd,enumerate,accents,pgffor,cite,wrapfig,multicol,float}
\usepackage{bibspacing}

\allowdisplaybreaks

\newcommand{\nc}{\newcommand}
\nc{\rnc}{\renewcommand}

\nc{\RP}{\mathcal{RP}}
\nc{\supp}{\operatorname{supp}}

\nc{\OEIS}{}

\usepackage{geometry}  \rnc{\ss}{\smallskip} \nc{\ms}{\medskip}  \nc{\nss}{\vspace{-3mm}} 

\makeatletter
\newenvironment{myalign}{%
  \def\align@preamble{
    &\hfil
     \strut@
     \setboxz@h{\@lign$\m@th\displaystyle{####}$}%
     \ifmeasuring@\savefieldlength@\fi
     \set@field
     \hfil
     \tabskip\z@skip
    &\setboxz@h{\@lign$\m@th\displaystyle{{}####}$}%
     \ifmeasuring@\savefieldlength@\fi
     \set@field
     \hfil
     \tabskip\alignsep@
  }\align}{\endalign}

\DeclareMathSymbol{\widehatsym}{\mathord}{largesymbols}{"62}
\newcommand\lowerwidehatsym{%
  \text{\smash{\raisebox{-1.3ex}{%
    $\widehatsym$}}}}
\newcommand\fixwidehat[1]{%
  \mathchoice
    {\accentset{\displaystyle\lowerwidehatsym}{#1}}
    {\accentset{\textstyle\lowerwidehatsym}{#1}}
    {\accentset{\scriptstyle\lowerwidehatsym}{#1}}
    {\accentset{\scriptscriptstyle\lowerwidehatsym}{#1}}
}
\rnc{\widehat}{\fixwidehat}

\begin{document}


\nc{\PnSn}{\P_n\sm\S_n}
\nc{\RPnSn}{\RP_n\sm\S_n}
\nc{\bbR}{\mathbb R}
\nc{\bbS}{\mathbb S}
\nc{\bbA}{\mathbb A}
\nc{\bbB}{\mathbb B}
\nc{\pre}{\preceq}
\nc{\Sub}{\operatorname{Sub}}
\nc{\MSA}{M(\S_n,\bbA_n)}
\nc{\Th}{\Theta}

\nc{\uvertcol}[2]{\fill[#2] (#1,2)circle(.2);}
\nc{\lvertcol}[2]{\fill[#2] (#1,0)circle(.2);}
\nc{\uvertcols}[2]{\foreach \x in {#1}{ \uvertcol{\x}{#2}}}
\nc{\lvertcols}[2]{\foreach \x in {#1}{ \lvertcol{\x}{#2}}}
\nc{\uverts}[1]{\foreach \x in {#1}{ \uvert{\x}}}
\nc{\lverts}[1]{\foreach \x in {#1}{ \lvert{\x}}}
\nc{\uvws}[1]{\foreach \x in {#1}{ \uvw{\x}}}
\nc{\lvws}[1]{\foreach \x in {#1}{ \lvw{\x}}}

\nc{\PBnt}{R^\tau[\PB_n]}
\nc{\Mnt}{R^\tau[\M_n]}
\nc{\PBntC}{\bbC^\tau[\PB_n]}
\nc{\MntC}{\bbC^\tau[\M_n]}
\nc{\bbC}{\mathbb C}
\nc{\RSr}{R[\S_{\br}]}
\nc{\fs}{\mathfrak s}
\nc{\ft}{\mathfrak t}
\nc{\fu}{\mathfrak u}
\nc{\fv}{\mathfrak v}
\nc{\rad}{\operatorname{rad}}
\nc{\dominates}{\unrhd}

\nc{\ubluebox}[2]{\bluebox{#1}{1.7}{#2}2\udotted{#1}{#2}}
\nc{\lbluebox}[2]{\bluebox{#1}0{#2}{.3}\ldotted{#1}{#2}}
\nc{\ublueboxes}[1]{{
\foreach \x/\y in {#1}
{ \ubluebox{\x}{\y}}}
}
\nc{\lblueboxes}[1]{{
\foreach \x/\y in {#1}
{ \lbluebox{\x}{\y}}}
}

\nc{\bluebox}[4]{
\draw[color=blue!20, fill=blue!20] (#1,#2)--(#3,#2)--(#3,#4)--(#1,#4)--(#1,#2);
}
\nc{\redbox}[4]{
\draw[color=red!20, fill=red!20] (#1,#2)--(#3,#2)--(#3,#4)--(#1,#4)--(#1,#2);
}

\nc{\bluetrap}[8]{
\draw[color=blue!20, fill=blue!20] (#1,#2)--(#3,#4)--(#5,#6)--(#7,#8)--(#1,#2);
}
\nc{\redtrap}[8]{
\draw[color=red!20, fill=red!20] (#1,#2)--(#3,#4)--(#5,#6)--(#7,#8)--(#1,#2);
}

\usetikzlibrary{decorations.markings}
\usetikzlibrary{arrows,matrix}
\usepgflibrary{arrows}
\tikzset{->-/.style={decoration={
  markings,
  mark=at position #1 with {\arrow{>}}},postaction={decorate}}}
\tikzset{-<-/.style={decoration={
  markings,
  mark=at position #1 with {\arrow{<}}},postaction={decorate}}}
\nc{\Unode}[1]{\draw(#1,-2)node{$U$};}
\nc{\Dnode}[1]{\draw(#1,-2)node{$D$};}
\nc{\Fnode}[1]{\draw(#1,-2)node{$F$};}
\nc{\Cnode}[1]{\draw(#1-.1,-2)node{$\phantom{+0},$};}
\nc{\Unodes}[1]{\foreach \x in {#1}{ \Unode{\x} }}
\nc{\Dnodes}[1]{\foreach \x in {#1}{ \Dnode{\x} }}
\nc{\Fnodes}[1]{\foreach \x in {#1}{ \Fnode{\x} }}
\nc{\Cnodes}[1]{\foreach \x in {#1}{ \Cnode{\x} }}
\nc{\Uedge}[2]{\draw[->-=0.6,line width=.3mm](#1,#2-9)--(#1+1,#2+1-9); \vertsm{#1}{#2-9} \vertsm{#1+1}{#2+1-9}}
\nc{\Dedge}[2]{\draw[->-=0.6,line width=.3mm](#1,#2-9)--(#1+1,#2-1-9); \vertsm{#1}{#2-9} \vertsm{#1+1}{#2-1-9}}
\nc{\Fedge}[2]{\draw[->-=0.6,line width=.3mm](#1,#2-9)--(#1+1,#2-9); \vertsm{#1}{#2-9} \vertsm{#1+1}{#2-9}}
\nc{\Uedges}[1]{\foreach \x/\y in {#1}{\Uedge{\x}{\y}}}
\nc{\Dedges}[1]{\foreach \x/\y in {#1}{\Dedge{\x}{\y}}}
\nc{\Fedges}[1]{\foreach \x/\y in {#1}{\Fedge{\x}{\y}}}
\nc{\xvertlabel}[1]{\draw(#1,-10+.6)node{{\tiny $#1$}};}
\nc{\yvertlabel}[1]{\draw(0-.4,-9+#1)node{{\tiny $#1$}};}
\nc{\xvertlabels}[1]{\foreach \x in {#1}{ \xvertlabel{\x} }}
\nc{\yvertlabels}[1]{\foreach \x in {#1}{ \yvertlabel{\x} }}

\nc{\bbE}{\mathbb E}
\nc{\floorn}{\lfloor\tfrac n2\rfloor}
\rnc{\sp}{\supseteq}
\rnc{\arraystretch}{1.2}

\nc{\bn}{{[n]}} \nc{\bt}{{[t]}} \nc{\ba}{{[a]}} \nc{\bl}{{[l]}} \nc{\bm}{{[m]}} \nc{\bk}{{[k]}} \nc{\br}{{[r]}} \nc{\bs}{{[s]}} \nc{\bnf}{{[n-1]}}

\nc{\M}{\mathcal M}
\nc{\G}{\mathcal G}
\nc{\F}{\mathfrak F}
\nc{\MnJ}{\mathcal M_n^J}
\nc{\EnJ}{\mathcal E_n^J}
\nc{\Mat}{\operatorname{Mat}}
\nc{\RegMnJ}{\Reg(\MnJ)}
\nc{\row}{\mathfrak r}
\nc{\col}{\mathfrak c}
\nc{\Row}{\operatorname{Row}}
\nc{\Col}{\operatorname{Col}}
\nc{\Span}{\operatorname{span}}
\nc{\mat}[4]{\left[\begin{matrix}#1&#2\\#3&#4\end{matrix}\right]}
\nc{\tmat}[4]{\left[\begin{smallmatrix}#1&#2\\#3&#4\end{smallmatrix}\right]}
\nc{\ttmat}[4]{{\tiny \left[\begin{smallmatrix}#1&#2\\#3&#4\end{smallmatrix}\right]}}
\nc{\tmatt}[9]{\left[\begin{smallmatrix}#1&#2&#3\\#4&#5&#6\\#7&#8&#9\end{smallmatrix}\right]}
\nc{\ttmatt}[9]{{\tiny \left[\begin{smallmatrix}#1&#2&#3\\#4&#5&#6\\#7&#8&#9\end{smallmatrix}\right]}}
\nc{\MnGn}{\M_n\sm\G_n}
\nc{\MrGr}{\M_r\sm\G_r}
\nc{\qbin}[2]{\left[\begin{matrix}#1\\#2\end{matrix}\right]_q}
\nc{\tqbin}[2]{\left[\begin{smallmatrix}#1\\#2\end{smallmatrix}\right]_q}
\nc{\qbinx}[3]{\left[\begin{matrix}#1\\#2\end{matrix}\right]_{#3}}
\nc{\tqbinx}[3]{\left[\begin{smallmatrix}#1\\#2\end{smallmatrix}\right]_{#3}}
\nc{\MNJ}{\M_nJ}
\nc{\JMN}{J\M_n}
\nc{\RegMNJ}{\Reg(\MNJ)}
\nc{\RegJMN}{\Reg(\JMN)}
\nc{\RegMMNJ}{\Reg(\MMNJ)}
\nc{\RegJMMN}{\Reg(\JMMN)}
\nc{\Wb}{\overline{W}}
\nc{\Xb}{\overline{X}}
\nc{\Yb}{\overline{Y}}
\nc{\Zb}{\overline{Z}}
\nc{\Sib}{\overline{\Si}}
\nc{\Om}{\Omega}
\nc{\Omb}{\overline{\Om}}
\nc{\Gab}{\overline{\Ga}}
\nc{\qfact}[1]{[#1]_q!}
\nc{\smat}[2]{\left[\begin{matrix}#1&#2\end{matrix}\right]}
\nc{\tsmat}[2]{\left[\begin{smallmatrix}#1&#2\end{smallmatrix}\right]}
\nc{\hmat}[2]{\left[\begin{matrix}#1\\#2\end{matrix}\right]}
\nc{\thmat}[2]{\left[\begin{smallmatrix}#1\\#2\end{smallmatrix}\right]}
\nc{\LVW}{\mathcal L(V,W)}
\nc{\KVW}{\mathcal K(V,W)}
\nc{\LV}{\mathcal L(V)}
\nc{\RegLVW}{\Reg(\LVW)}
\nc{\sM}{\mathscr M}
\nc{\sN}{\mathscr N}
\rnc{\iff}{\ \Leftrightarrow\ }
\nc{\Hom}{\operatorname{Hom}}
\nc{\End}{\operatorname{End}}
\nc{\Aut}{\operatorname{Aut}}
\nc{\Lin}{\mathcal L}
\nc{\Hommn}{\Hom(V_m,V_n)}
\nc{\Homnm}{\Hom(V_n,V_m)}
\nc{\Homnl}{\Hom(V_n,V_l)}
\nc{\Homkm}{\Hom(V_k,V_m)}
\nc{\Endm}{\End(V_m)}
\nc{\Endn}{\End(V_n)}
\nc{\Endr}{\End(V_r)}
\nc{\Autm}{\Aut(V_m)}
\nc{\Autn}{\Aut(V_n)}
\nc{\MmnJ}{\M_{mn}^J}
\nc{\MmnA}{\M_{mn}^A}
\nc{\MmnB}{\M_{mn}^B}
\nc{\Mmn}{\M_{mn}}
\nc{\Mkl}{\M_{kl}}
\nc{\Mnm}{\M_{nm}}
\nc{\EmnJ}{\mathcal E_{mn}^J}
\nc{\MmGm}{\M_m\sm\G_m}
\nc{\RegMmnJ}{\Reg(\MmnJ)}
\rnc{\implies}{\ \Rightarrow\ }
\nc{\DMmn}[1]{D_{#1}(\Mmn)}
\nc{\DMmnJ}[1]{D_{#1}(\MmnJ)}
\nc{\MMNJ}{\Mmn J}
\nc{\JMMN}{J\Mmn}
\nc{\JMMNJ}{J\Mmn J}
\nc{\Inr}{\mathcal I(V_n,W_r)}
\nc{\Lnr}{\mathcal L(V_n,W_r)}
\nc{\Knr}{\mathcal K(V_n,W_r)}
\nc{\Imr}{\mathcal I(V_m,W_r)}
\nc{\Kmr}{\mathcal K(V_m,W_r)}
\nc{\Lmr}{\mathcal L(V_m,W_r)}
\nc{\Kmmr}{\mathcal K(V_m,W_{m-r})}
\nc{\tr}{{\operatorname{T}}}
\nc{\MMN}{\MmnA(\F_1)}
\nc{\MKL}{\Mkl^B(\F_2)}
\nc{\RegMMN}{\Reg(\MmnA(\F_1))}
\nc{\RegMKL}{\Reg(\Mkl^B(\F_2))}
\nc{\gRhA}{\widehat{\mathscr R}^A}
\nc{\gRhB}{\widehat{\mathscr R}^B}
\nc{\gLhA}{\widehat{\mathscr L}^A}
\nc{\gLhB}{\widehat{\mathscr L}^B}
\nc{\timplies}{\Rightarrow}
\nc{\tiff}{\Leftrightarrow}
\nc{\Sija}{S_{ij}^a}
\nc{\dmat}[8]{\draw(#1*1.5,#2)node{$\left[\begin{smallmatrix}#3&#4&#5\\#6&#7&#8\end{smallmatrix}\right]$};}
\nc{\bdmat}[8]{\draw(#1*1.5,#2)node{${\mathbf{\left[\begin{smallmatrix}#3&#4&#5\\#6&#7&#8\end{smallmatrix}\right]}}$};}
\nc{\rdmat}[8]{\draw(#1*1.5,#2)node{\rotatebox{90}{$\left[\begin{smallmatrix}#3&#4&#5\\#6&#7&#8\end{smallmatrix}\right]$}};}
\nc{\rldmat}[8]{\draw(#1*1.5-0.375,#2)node{\rotatebox{90}{$\left[\begin{smallmatrix}#3&#4&#5\\#6&#7&#8\end{smallmatrix}\right]$}};}
\nc{\rrdmat}[8]{\draw(#1*1.5+.375,#2)node{\rotatebox{90}{$\left[\begin{smallmatrix}#3&#4&#5\\#6&#7&#8\end{smallmatrix}\right]$}};}
\nc{\rfldmat}[8]{\draw(#1*1.5-0.375+.15,#2)node{\rotatebox{90}{$\left[\begin{smallmatrix}#3&#4&#5\\#6&#7&#8\end{smallmatrix}\right]$}};}
\nc{\rfrdmat}[8]{\draw(#1*1.5+.375-.15,#2)node{\rotatebox{90}{$\left[\begin{smallmatrix}#3&#4&#5\\#6&#7&#8\end{smallmatrix}\right]$}};}
\nc{\xL}{[x]_{\! _\gL}}\nc{\yL}{[y]_{\! _\gL}}\nc{\xR}{[x]_{\! _\gR}}\nc{\yR}{[y]_{\! _\gR}}\nc{\xH}{[x]_{\! _\gH}}\nc{\yH}{[y]_{\! _\gH}}\nc{\XK}{[X]_{\! _\gK}}\nc{\xK}{[x]_{\! _\gK}}
\nc{\RegSija}{\Reg(\Sija)}
\nc{\MnmK}{\M_{nm}^K}
\nc{\cC}{\mathcal C}
\nc{\cR}{\mathcal R}
\nc{\Ckl}{\cC_k(l)}
\nc{\Rkl}{\cR_k(l)}
\nc{\Cmr}{\cC_m(r)}
\nc{\Rmr}{\cR_m(r)}
\nc{\Cnr}{\cC_n(r)}
\nc{\Rnr}{\cR_n(r)}
\nc{\Z}{\mathbb Z}

\nc{\Reg}{\operatorname{Reg}}
\nc{\TXa}{\T_X^a}
\nc{\TXA}{\T(X,A)}
\nc{\TXal}{\T(X,\al)}
\nc{\RegTXa}{\Reg(\TXa)}
\nc{\RegTXA}{\Reg(\TXA)}
\nc{\RegTXal}{\Reg(\TXal)}
\nc{\PalX}{\P_\al(X)}
\nc{\EAX}{\E_A(X)}
\nc{\Bb}{\overline{B}}
\nc{\bb}{\overline{\be}}
\nc{\bw}{{\bf w}}
\nc{\bz}{{\bf z}}
\nc{\TASA}{\T_A\sm\S_A}
\nc{\Ub}{\overline{U}}
\nc{\Vb}{\overline{V}}
\nc{\eb}{\overline{e}}
\nc{\ob}{\overline{o}}
\nc{\cb}{\overline{c}}
\nc{\db}{\overline{d}}
\nc{\tb}{\overline{t}}
\nc{\Eb}{\overline{E}}
\nc{\Ob}{\overline{O}}
\nc{\Qb}{\overline{Q}}
\nc{\EXa}{\E_X^a}
\nc{\oijr}{1\leq i<j\leq r}
\nc{\veb}{\overline{\ve}}
\nc{\bbT}{\mathbb T}
\nc{\Surj}{\operatorname{Surj}}
\nc{\Sone}{S^{(1)}}
\nc{\fillbox}[2]{\draw[fill=gray!30](#1,#2)--(#1+1,#2)--(#1+1,#2+1)--(#1,#2+1)--(#1,#2);}
\nc{\raa}{\rangle_J}
\nc{\raJ}{\rangle_J}
\nc{\Ea}{E_J}
\nc{\EJ}{E_J}
\nc{\ep}{\epsilon} \nc{\ve}{\varepsilon}
\nc{\IXa}{\I_X^a}
\nc{\RegIXa}{\Reg(\IXa)}
\nc{\JXa}{\J_X^a}
\nc{\RegJXa}{\Reg(\JXa)}
\nc{\IXA}{\I(X,A)}
\nc{\IAX}{\I(A,X)}
\nc{\RegIXA}{\Reg(\IXA)}
\nc{\RegIAX}{\Reg(\IAX)}
\nc{\trans}[2]{\left(\begin{smallmatrix} #1 \\ #2 \end{smallmatrix}\right)}
\nc{\bigtrans}[2]{\left(\begin{matrix} #1 \\ #2 \end{matrix}\right)}
\nc{\lmap}[1]{\mapstochar \xrightarrow {\ #1\ }}
\nc{\EaTXa}{E}

\nc{\gL}{\mathscr L}
\nc{\gR}{\mathscr R}
\nc{\gH}{\mathscr H}
\nc{\gJ}{\mathscr J}
\nc{\gD}{\mathscr D}
\nc{\gK}{\mathscr K}
\nc{\gLa}{\mathscr L^a}
\nc{\gRa}{\mathscr R^a}
\nc{\gHa}{\mathscr H^a}
\nc{\gJa}{\mathscr J^a}
\nc{\gDa}{\mathscr D^a}
\nc{\gKa}{\mathscr K^a}
\nc{\gLJ}{\mathscr L^J}
\nc{\gRJ}{\mathscr R^J}
\nc{\gHJ}{\mathscr H^J}
\nc{\gJJ}{\mathscr J^J}
\nc{\gDJ}{\mathscr D^J}
\nc{\gKJ}{\mathscr K^J}
\nc{\gLh}{\widehat{\mathscr L}^J}
\nc{\gRh}{\widehat{\mathscr R}^J}
\nc{\gHh}{\widehat{\mathscr H}^J}
\nc{\gJh}{\widehat{\mathscr J}^J}
\nc{\gDh}{\widehat{\mathscr D}^J}
\nc{\gKh}{\widehat{\mathscr K}^J}
\nc{\Lh}{\widehat{L}^J}
\nc{\Rh}{\widehat{R}}
\nc{\Hh}{\widehat{H}^J}
\nc{\Jh}{\widehat{J}^J}
\nc{\Dh}{\widehat{D}^J}
\nc{\Kh}{\widehat{K}^J}
\nc{\gLb}{\widehat{\mathscr L}}
\nc{\gRb}{\widehat{\mathscr R}}
\nc{\gHb}{\widehat{\mathscr H}}
\nc{\gJb}{\widehat{\mathscr J}}
\nc{\gDb}{\widehat{\mathscr D}}
\nc{\gKb}{\widehat{\mathscr K}}
\nc{\Lb}{\widehat{L}^J}
\nc{\Rb}{\widehat{R}^J}
\nc{\Hb}{\widehat{H}^J}
\nc{\Jb}{\widehat{J}^J}
\nc{\Db}{\overline{D}}
\nc{\Kb}{\widehat{K}}

\hyphenation{mon-oid mon-oids}

\nc{\itemit}[1]{\item[\emph{(#1)}]}
\nc{\E}{\mathcal E}
\nc{\TX}{\T(X)}
\nc{\TXP}{\T(X,\P)}
\nc{\EX}{\E(X)}
\nc{\EXP}{\E(X,\P)}
\nc{\SX}{\S(X)}
\nc{\SXP}{\S(X,\P)}
\nc{\Sing}{\operatorname{Sing}}
\nc{\idrank}{\operatorname{idrank}}
\nc{\SingXP}{\Sing(X,\P)}
\nc{\De}{\Delta}
\nc{\sgp}{\operatorname{sgp}}
\nc{\mon}{\operatorname{mon}}
\nc{\Dn}{\mathcal D_n}
\nc{\Dm}{\mathcal D_m}

\nc{\lline}[1]{\draw(3*#1,0)--(3*#1+2,0);}
\nc{\uline}[1]{\draw(3*#1,5)--(3*#1+2,5);}
\nc{\thickline}[2]{\draw(3*#1,5)--(3*#2,0); \draw(3*#1+2,5)--(3*#2+2,0) ;}
\nc{\thicklabel}[3]{\draw(3*#1+1+3*#2*0.15-3*#1*0.15,4.25)node{{\tiny $#3$}};}

\nc{\slline}[3]{\draw(3*#1+#3,0+#2)--(3*#1+2+#3,0+#2);}
\nc{\suline}[3]{\draw(3*#1+#3,5+#2)--(3*#1+2+#3,5+#2);}
\nc{\sthickline}[4]{\draw(3*#1+#4,5+#3)--(3*#2+#4,0+#3); \draw(3*#1+2+#4,5+#3)--(3*#2+2+#4,0+#3) ;}
\nc{\sthicklabel}[5]{\draw(3*#1+1+3*#2*0.15-3*#1*0.15+#5,4.25+#4)node{{\tiny $#3$}};}

\nc{\stll}[5]{\sthickline{#1}{#2}{#4}{#5} \sthicklabel{#1}{#2}{#3}{#4}{#5}}
\nc{\tll}[3]{\stll{#1}{#2}{#3}00}

\nc{\mfourpic}[9]{
\slline1{#9}0
\slline3{#9}0
\slline4{#9}0
\slline5{#9}0
\suline1{#9}0
\suline3{#9}0
\suline4{#9}0
\suline5{#9}0
\stll1{#1}{#5}{#9}{0}
\stll3{#2}{#6}{#9}{0}
\stll4{#3}{#7}{#9}{0}
\stll5{#4}{#8}{#9}{0}
\draw[dotted](6,0+#9)--(8,0+#9);
\draw[dotted](6,5+#9)--(8,5+#9);
}
\nc{\vdotted}[1]{
\draw[dotted](3*#1,10)--(3*#1,15);
\draw[dotted](3*#1+2,10)--(3*#1+2,15);
}

\nc{\Clab}[2]{
\sthicklabel{#1}{#1}{{}_{\phantom{#1}}C_{#1}}{1.25+5*#2}0
}
\nc{\sClab}[3]{
\sthicklabel{#1}{#1}{{}_{\phantom{#1}}C_{#1}}{1.25+5*#2}{#3}
}
\nc{\Clabl}[3]{
\sthicklabel{#1}{#1}{{}_{\phantom{#3}}C_{#3}}{1.25+5*#2}0
}
\nc{\sClabl}[4]{
\sthicklabel{#1}{#1}{{}_{\phantom{#4}}C_{#4}}{1.25+5*#2}{#3}
}
\nc{\Clabll}[3]{
\sthicklabel{#1}{#1}{C_{#3}}{1.25+5*#2}0
}
\nc{\sClabll}[4]{
\sthicklabel{#1}{#1}{C_{#3}}{1.25+5*#2}{#3}
}

\nc{\mtwopic}[6]{
\slline1{#6*5}{#5}
\slline2{#6*5}{#5}
\suline1{#6*5}{#5}
\suline2{#6*5}{#5}
\stll1{#1}{#3}{#6*5}{#5}
\stll2{#2}{#4}{#6*5}{#5}
}
\nc{\mtwopicl}[6]{
\slline1{#6*5}{#5}
\slline2{#6*5}{#5}
\suline1{#6*5}{#5}
\suline2{#6*5}{#5}
\stll1{#1}{#3}{#6*5}{#5}
\stll2{#2}{#4}{#6*5}{#5}
\sClabl1{#6}{#5}{i}
\sClabl2{#6}{#5}{j}
}

\nc{\keru}{\operatorname{ker}^\wedge} \nc{\kerl}{\operatorname{ker}_\vee}

\nc{\coker}{\operatorname{coker}}
\nc{\KER}{\ker}
\nc{\N}{\mathbb N}
\nc{\LaBn}{L_\al(\B_n)}
\nc{\RaBn}{R_\al(\B_n)}
\nc{\LaPBn}{L_\al(\PB_n)}
\nc{\RaPBn}{R_\al(\PB_n)}
\nc{\rhorBn}{\rho_r(\B_n)}
\nc{\DrBn}{D_r(\B_n)}
\nc{\DrPn}{D_r(\P_n)}
\nc{\DrPBn}{D_r(\PB_n)}
\nc{\DrKn}{D_r(\K_n)}
\nc{\alb}{\al_{\vee}}
\nc{\beb}{\be^{\wedge}}
\nc{\Bal}{\operatorname{Bal}}
\nc{\Red}{\operatorname{Red}}
\nc{\Pnxi}{\P_n^\xi}
\nc{\Bnxi}{\B_n^\xi}
\nc{\PBnxi}{\PB_n^\xi}
\nc{\Knxi}{\K_n^\xi}
\nc{\C}{\mathbb C}
\nc{\exi}{e^\xi}
\nc{\Exi}{E^\xi}
\nc{\eximu}{e^\xi_\mu}
\nc{\Eximu}{E^\xi_\mu}
\nc{\REF}{ {\red [Ref?]} }
\nc{\GL}{\operatorname{GL}}
\rnc{\O}{\mathcal O}

\nc{\vtx}[2]{\fill (#1,#2)circle(.2);}
\nc{\lvtx}[2]{\fill (#1,0)circle(.2);}
\nc{\uvtx}[2]{\fill (#1,1.5)circle(.2);}

\nc{\Eq}{\mathfrak{Eq}}
\nc{\Gau}{\Ga^\wedge} \nc{\Gal}{\Ga_\vee}
\nc{\Lamu}{\Lam^\wedge} \nc{\Laml}{\Lam_\vee}
\nc{\bX}{{\bf X}}
\nc{\bY}{{\bf Y}}
\nc{\ds}{\displaystyle}

\nc{\uuvert}[1]{\fill (#1,3)circle(.2);}
\nc{\uuuvert}[1]{\fill (#1,4.5)circle(.2);}
\nc{\overt}[1]{\fill (#1,0)circle(.1);}
\nc{\overtl}[3]{\node[vertex] (#3) at (#1,0) {  {\tiny $#2$} };}
\nc{\cv}[2]{\draw(#1,1.5) to [out=270,in=90] (#2,0);}
\nc{\cvs}[2]{\draw(#1,1.5) to [out=270+30,in=90+30] (#2,0);}
\nc{\ucv}[2]{\draw(#1,3) to [out=270,in=90] (#2,1.5);}
\nc{\uucv}[2]{\draw(#1,4.5) to [out=270,in=90] (#2,3);}
\nc{\textpartn}[1]{{\lower1.0 ex\hbox{\begin{tikzpicture}[xscale=.3,yscale=0.3] #1 \end{tikzpicture}}}}
\nc{\textpartnx}[2]{{\lower1.0 ex\hbox{\begin{tikzpicture}[xscale=.3,yscale=0.3] 
\foreach \x in {1,...,#1}
{ \uvert{\x} \lvert{\x} }
#2 \end{tikzpicture}}}}
\nc{\disppartnx}[2]{{\lower1.0 ex\hbox{\begin{tikzpicture}[scale=0.3] 
\foreach \x in {1,...,#1}
{ \uvert{\x} \lvert{\x} }
#2 \end{tikzpicture}}}}
\nc{\disppartnxd}[2]{{\lower2.1 ex\hbox{\begin{tikzpicture}[scale=0.3] 
\foreach \x in {1,...,#1}
{ \uuvert{\x} \uvert{\x} \lvert{\x} }
#2 \end{tikzpicture}}}}
\nc{\disppartnxdn}[2]{{\lower2.1 ex\hbox{\begin{tikzpicture}[scale=0.3] 
\foreach \x in {1,...,#1}
{ \uuvert{\x} \lvert{\x} }
#2 \end{tikzpicture}}}}
\nc{\disppartnxdd}[2]{{\lower3.6 ex\hbox{\begin{tikzpicture}[scale=0.3] 
\foreach \x in {1,...,#1}
{ \uuuvert{\x} \uuvert{\x} \uvert{\x} \lvert{\x} }
#2 \end{tikzpicture}}}}

\nc{\dispgax}[2]{{\lower0.0 ex\hbox{\begin{tikzpicture}[scale=0.3] 
#2
\foreach \x in {1,...,#1}
{\lvert{\x} }
 \end{tikzpicture}}}}
\nc{\textgax}[2]{{\lower0.4 ex\hbox{\begin{tikzpicture}[scale=0.3] 
#2
\foreach \x in {1,...,#1}
{\lvert{\x} }
 \end{tikzpicture}}}}
\nc{\textlinegraph}[2]{{\raise#1 ex\hbox{\begin{tikzpicture}[scale=0.8] 
#2
 \end{tikzpicture}}}}
\nc{\textlinegraphl}[2]{{\raise#1 ex\hbox{\begin{tikzpicture}[scale=0.8] 
\tikzstyle{vertex}=[circle,draw=black, fill=white, inner sep = 0.07cm]
#2
 \end{tikzpicture}}}}
\nc{\displinegraph}[1]{{\lower0.0 ex\hbox{\begin{tikzpicture}[scale=0.6] 
#1
 \end{tikzpicture}}}}
 
\nc{\disppartnthreeone}[1]{{\lower1.0 ex\hbox{\begin{tikzpicture}[scale=0.3] 
\foreach \x in {1,2,3,5,6}
{ \uvert{\x} }
\foreach \x in {1,2,4,5,6}
{ \lvert{\x} }
\draw[dotted] (3.5,1.5)--(4.5,1.5);
\draw[dotted] (2.5,0)--(3.5,0);
#1 \end{tikzpicture}}}}

\nc{\partn}[4]{\left( \begin{array}{c|c} 
#1 \ & \ #3 \ \ \\ \cline{2-2}
#2 \ & \ #4 \ \
\end{array} \!\!\! \right)}
\nc{\partnlong}[6]{\partn{#1}{#2}{#3,\ #4}{#5,\ #6}} 
\nc{\partnsh}[2]{\left( \begin{array}{c} 
#1 \\
#2 
\end{array} \right)}
\nc{\partncodefz}[3]{\partn{#1}{#2}{#3}{\emptyset}}
\nc{\partndefz}[3]{{\partn{#1}{#2}{\emptyset}{#3}}}
\nc{\partnlast}[2]{\left( \begin{array}{c|c}
#1 \ &  \ #2 \\
#1 \ &  \ #2
\end{array} \right)}

\nc{\partnlist}[8]{
\left( \begin{array}{c|c|c|c|c|c} 
\!\! #1 & \cdots & #2 & #3 & \cdots & #4\ \ \\ \cline{4-6}
\!\! #5 & \cdots & #6 & #7 & \cdots & #8 \ \
\end{array} \!\!\! \right)
}

\nc{\partnlistnodef}[8]{
\left( \begin{array}{c|c|c|c|c|c} 
\!\! #1 & \cdots & #2 & #3 & \cdots & #4\ \ \\ 
\!\! #5 & \cdots & #6 & #7 & \cdots & #8 \ \
\end{array} \!\!\! \right)
}

\nc{\partnlistnodefshort}[4]{
\left( \begin{array}{c|c|c} 
\!\! #1 & \cdots & #2 \\ 
\!\! #3 & \cdots & #4 
\end{array} \!\!\! \right)
}

\nc{\partnlistnodefshortish}[6]{
\left( \begin{array}{c|c|c|c} 
\!\! #1 & \cdots & #2 & #3 \\ 
\!\! #4 & \cdots & #5 & #6 
\end{array} \!\!\! \right)
}

\nc{\uuarcx}[3]{\draw(#1,3)arc(180:270:#3) (#1+#3,3-#3)--(#2-#3,3-#3) (#2-#3,3-#3) arc(270:360:#3);}
\nc{\uuarc}[2]{\uuarcx{#1}{#2}{.4}}
\nc{\uuuarcx}[3]{\draw(#1,4.5)arc(180:270:#3) (#1+#3,4.5-#3)--(#2-#3,4.5-#3) (#2-#3,4.5-#3) arc(270:360:#3);}
\nc{\uuuarc}[2]{\uuuarcx{#1}{#2}{.4}}
\nc{\udarcx}[3]{\draw(#1,1.5)arc(180:90:#3) (#1+#3,1.5+#3)--(#2-#3,1.5+#3) (#2-#3,1.5+#3) arc(90:0:#3);}
\nc{\udarc}[2]{\udarcx{#1}{#2}{.4}}
\nc{\uudarcx}[3]{\draw(#1,3)arc(180:90:#3) (#1+#3,3+#3)--(#2-#3,3+#3) (#2-#3,3+#3) arc(90:0:#3);}
\nc{\uudarc}[2]{\uudarcx{#1}{#2}{.4}}
\nc{\darcxhalf}[3]{\draw(#1,0)arc(180:90:#3) (#1+#3,#3)--(#2,#3) ;}
\nc{\darchalf}[2]{\darcxhalf{#1}{#2}{.4}}
\nc{\uarcxhalf}[3]{\draw(#1,1.5)arc(180:270:#3) (#1+#3,1.5-#3)--(#2,1.5-#3) ;}
\nc{\uarchalf}[2]{\uarcxhalf{#1}{#2}{.4}}
\nc{\uarcxhalfr}[3]{\draw (#1+#3,1.5-#3)--(#2-#3,1.5-#3) (#2-#3,1.5-#3) arc(270:360:#3);}
\nc{\uarchalfr}[2]{\uarcxhalfr{#1}{#2}{.4}}

\nc{\bdarcx}[3]{\draw[blue](#1,0)arc(180:90:#3) (#1+#3,#3)--(#2-#3,#3) (#2-#3,#3) arc(90:0:#3);}
\nc{\bdarc}[2]{\darcx{#1}{#2}{.4}}
\nc{\rduarcx}[3]{\draw[red](#1,0)arc(180:270:#3) (#1+#3,0-#3)--(#2-#3,0-#3) (#2-#3,0-#3) arc(270:360:#3);}
\nc{\rduarc}[2]{\uarcx{#1}{#2}{.4}}
\nc{\duarcx}[3]{\draw(#1,0)arc(180:270:#3) (#1+#3,0-#3)--(#2-#3,0-#3) (#2-#3,0-#3) arc(270:360:#3);}
\nc{\duarc}[2]{\uarcx{#1}{#2}{.4}}

\nc{\uuv}[1]{\fill (#1,4)circle(.1);}
\nc{\uv}[1]{\fill (#1,2)circle(.1);}
\nc{\lv}[1]{\fill (#1,0)circle(.1);}
\nc{\uvw}[1]{\draw[fill=white] (#1,2)circle(.18);}
\nc{\lvw}[1]{\draw (#1,0)circle(.18);}
\nc{\uvred}[1]{\fill[red] (#1,2)circle(.1);}
\nc{\lvred}[1]{\fill[red] (#1,0)circle(.1);}
\nc{\lvwhite}[1]{\fill[white] (#1,0)circle(.1);}
\nc{\buv}[1]{\fill (#1,2)circle(.18);}
\nc{\blv}[1]{\fill (#1,0)circle(.18);}

\nc{\uvs}[1]{{
\foreach \x in {#1}
{ \uv{\x}}
}}
\nc{\uuvs}[1]{{
\foreach \x in {#1}
{ \uuv{\x}}
}}
\nc{\lvs}[1]{{
\foreach \x in {#1}
{ \lv{\x}}
}}

\nc{\buvs}[1]{{
\foreach \x in {#1}
{ \buv{\x}}
}}
\nc{\blvs}[1]{{
\foreach \x in {#1}
{ \blv{\x}}
}}

\nc{\uvreds}[1]{{
\foreach \x in {#1}
{ \uvred{\x}}
}}
\nc{\lvreds}[1]{{
\foreach \x in {#1}
{ \lvred{\x}}
}}

\nc{\uudotted}[2]{\draw [dotted] (#1,4)--(#2,4);}
\nc{\uudotteds}[1]{{
\foreach \x/\y in {#1}
{ \uudotted{\x}{\y}}
}}
\nc{\uudottedsm}[2]{\draw [dotted] (#1+.4,4)--(#2-.4,4);}
\nc{\uudottedsms}[1]{{
\foreach \x/\y in {#1}
{ \uudottedsm{\x}{\y}}
}}
\nc{\udottedsm}[2]{\draw [dotted] (#1+.4,2)--(#2-.4,2);}
\nc{\udottedsms}[1]{{
\foreach \x/\y in {#1}
{ \udottedsm{\x}{\y}}
}}
\nc{\udotted}[2]{\draw [dotted] (#1,2)--(#2,2);}
\nc{\udotteds}[1]{{
\foreach \x/\y in {#1}
{ \udotted{\x}{\y}}
}}
\nc{\ldotted}[2]{\draw [dotted] (#1,0)--(#2,0);}
\nc{\ldotteds}[1]{{
\foreach \x/\y in {#1}
{ \ldotted{\x}{\y}}
}}
\nc{\ldottedsm}[2]{\draw [dotted] (#1+.4,0)--(#2-.4,0);}
\nc{\ldottedsms}[1]{{
\foreach \x/\y in {#1}
{ \ldottedsm{\x}{\y}}
}}

\nc{\stlinest}[2]{\draw(#1,4)--(#2,0);}

\nc{\stlined}[2]{\draw[dotted](#1,2)--(#2,0);}

\nc{\tlab}[2]{\draw(#1,2)node[above]{\tiny $#2$};}
\nc{\tudots}[1]{\draw(#1,2)node{$\cdots$};}
\nc{\tldots}[1]{\draw(#1,0)node{$\cdots$};}

\nc{\huv}[1]{\fill (#1,1)circle(.1);}
\nc{\llv}[1]{\fill (#1,-2)circle(.1);}
\nc{\arcup}[2]{
\draw(#1,2)arc(180:270:.4) (#1+.4,1.6)--(#2-.4,1.6) (#2-.4,1.6) arc(270:360:.4);
}
\nc{\harcup}[2]{
\draw(#1,1)arc(180:270:.4) (#1+.4,.6)--(#2-.4,.6) (#2-.4,.6) arc(270:360:.4);
}
\nc{\arcdn}[2]{
\draw(#1,0)arc(180:90:.4) (#1+.4,.4)--(#2-.4,.4) (#2-.4,.4) arc(90:0:.4);
}
\nc{\cve}[2]{
\draw(#1,2) to [out=270,in=90] (#2,0);
}
\nc{\hcve}[2]{
\draw(#1,1) to [out=270,in=90] (#2,0);
}
\nc{\catarc}[3]{
\draw(#1,2)arc(180:270:#3) (#1+#3,2-#3)--(#2-#3,2-#3) (#2-#3,2-#3) arc(270:360:#3);
}

\nc{\arcr}[2]{
\draw[red](#1,0)arc(180:90:.4) (#1+.4,.4)--(#2-.4,.4) (#2-.4,.4) arc(90:0:.4);
}
\nc{\arcb}[2]{
\draw[blue](#1,2-2)arc(180:270:.4) (#1+.4,1.6-2)--(#2-.4,1.6-2) (#2-.4,1.6-2) arc(270:360:.4);
}
\nc{\loopr}[1]{
\draw[blue](#1,-2) edge [out=130,in=50,loop] ();
}
\nc{\loopb}[1]{
\draw[red](#1,-2) edge [out=180+130,in=180+50,loop] ();
}
\nc{\redto}[2]{\draw[red](#1,0)--(#2,0);}
\nc{\bluto}[2]{\draw[blue](#1,0)--(#2,0);}
\nc{\dotto}[2]{\draw[dotted](#1,0)--(#2,0);}
\nc{\lloopr}[2]{\draw[red](#1,0)arc(0:360:#2);}
\nc{\lloopb}[2]{\draw[blue](#1,0)arc(0:360:#2);}
\nc{\rloopr}[2]{\draw[red](#1,0)arc(-180:180:#2);}
\nc{\rloopb}[2]{\draw[blue](#1,0)arc(-180:180:#2);}
\nc{\uloopr}[2]{\draw[red](#1,0)arc(-270:270:#2);}
\nc{\uloopb}[2]{\draw[blue](#1,0)arc(-270:270:#2);}
\nc{\dloopr}[2]{\draw[red](#1,0)arc(-90:270:#2);}
\nc{\dloopb}[2]{\draw[blue](#1,0)arc(-90:270:#2);}
\nc{\llloopr}[2]{\draw[red](#1,0-2)arc(0:360:#2);}
\nc{\llloopb}[2]{\draw[blue](#1,0-2)arc(0:360:#2);}
\nc{\lrloopr}[2]{\draw[red](#1,0-2)arc(-180:180:#2);}
\nc{\lrloopb}[2]{\draw[blue](#1,0-2)arc(-180:180:#2);}
\nc{\ldloopr}[2]{\draw[red](#1,0-2)arc(-270:270:#2);}
\nc{\ldloopb}[2]{\draw[blue](#1,0-2)arc(-270:270:#2);}
\nc{\luloopr}[2]{\draw[red](#1,0-2)arc(-90:270:#2);}
\nc{\luloopb}[2]{\draw[blue](#1,0-2)arc(-90:270:#2);}

\nc{\larcb}[2]{
\draw[blue](#1,0-2)arc(180:90:.4) (#1+.4,.4-2)--(#2-.4,.4-2) (#2-.4,.4-2) arc(90:0:.4);
}
\nc{\larcr}[2]{
\draw[red](#1,2-2-2)arc(180:270:.4) (#1+.4,1.6-2-2)--(#2-.4,1.6-2-2) (#2-.4,1.6-2-2) arc(270:360:.4);
}

\rnc{\H}{\mathscr H}
\rnc{\L}{\mathscr L}
\nc{\R}{\mathcal R}
\nc{\D}{\mathscr D}
\nc{\J}{\mathcal J}

\nc{\ssim}{\mathrel{\raise0.25 ex\hbox{\oalign{$\approx$\crcr\noalign{\kern-0.84 ex}$\approx$}}}}
\nc{\POI}{\mathcal{O}}
\nc{\wb}{\overline{w}}
\nc{\ub}{\overline{u}}
\nc{\vb}{\overline{v}}
\nc{\fb}{\overline{f}}
\nc{\gb}{\overline{g}}
\nc{\hb}{\overline{h}}
\nc{\pb}{\overline{p}}
\nc{\xb}{\overline{x}}
\nc{\qb}{\overline{q}}
\rnc{\sb}{\overline{s}}
\nc{\Sb}{\overline{\si}}
\nc{\XR}{\pres{X}{R\,}}
\nc{\YQ}{\pres{Y}{Q}}
\nc{\ZP}{\pres{Z}{P\,}}
\nc{\XRone}{\pres{X_1}{R_1}}
\nc{\XRtwo}{\pres{X_2}{R_2}}
\nc{\XRthree}{\pres{X_1\cup X_2}{R_1\cup R_2\cup R_3}}
\nc{\er}{\eqref}
\nc{\larr}{\mathrel{\hspace{-0.35 ex}>\hspace{-1.1ex}-}\hspace{-0.35 ex}}
\nc{\rarr}{\mathrel{\hspace{-0.35 ex}-\hspace{-0.5ex}-\hspace{-2.3ex}>\hspace{-0.35 ex}}}
\nc{\lrarr}{\mathrel{\hspace{-0.35 ex}>\hspace{-1.1ex}-\hspace{-0.5ex}-\hspace{-2.3ex}>\hspace{-0.35 ex}}}
\nc{\nn}{\tag*{}}
\nc{\epfal}{\tag*{$\Box$}}
\nc{\tagd}[1]{\tag*{(#1)$'$}}
\nc{\ldb}{[\![}
\nc{\rdb}{]\!]}
\nc{\sm}{\setminus}
\nc{\I}{\mathcal I}
\nc{\InSn}{\I_n\setminus\S_n}
\nc{\dom}{\operatorname{dom}} \nc{\codom}{\operatorname{codom}}
\nc{\ojin}{1\leq j<i\leq n}
\nc{\eh}{\widehat{e}}
\nc{\wh}{\widehat{w}}
\nc{\uh}{\widehat{u}}
\nc{\vh}{\widehat{v}}
\nc{\sh}{\widehat{s}}
\nc{\fh}{\widehat{f}}
\nc{\textres}[1]{\text{\emph{#1}}}
\nc{\aand}{\emph{\ and \ }}
\nc{\iif}{\emph{\ if \ }}
\nc{\textlarr}{\ \larr\ }
\nc{\textrarr}{\ \rarr\ }
\nc{\textlrarr}{\ \lrarr\ }

\nc{\comma}{,\ }

\nc{\COMMA}{,\quad}
\nc{\TnSn}{\T_n\setminus\S_n} 
\nc{\TmSm}{\T_m\setminus\S_m} 
\nc{\TXSX}{\T_X\setminus\S_X} 
\rnc{\S}{\mathcal S}

\nc{\T}{\mathcal T} 
\nc{\A}{\mathscr A} 
\nc{\B}{\mathcal B} 
\rnc{\P}{\mathcal P} 
\nc{\K}{\mathcal K}
\nc{\PB}{\mathcal{PB}} 
\nc{\rank}{\operatorname{rank}}

\nc{\mtt}{\!\!\!\mt\!\!\!}

\nc{\sub}{\subseteq}
\nc{\la}{\langle}
\nc{\ra}{\rangle}
\nc{\mt}{\mapsto}
\nc{\im}{\mathrm{im}}
\nc{\id}{\mathrm{id}}
\nc{\al}{\alpha}
\nc{\be}{\beta}
\nc{\ga}{\gamma}
\nc{\Ga}{\Gamma}
\nc{\de}{\delta}
\nc{\ka}{\kappa}
\nc{\lam}{\lambda}
\nc{\Lam}{\Lambda}
\nc{\si}{\sigma}
\nc{\Si}{\Sigma}
\nc{\oijn}{1\leq i<j\leq n}
\nc{\oijm}{1\leq i<j\leq m}

\nc{\comm}{\rightleftharpoons}
\nc{\AND}{\qquad\text{and}\qquad}

\nc{\bit}{\vspace{-3 truemm}\begin{itemize}}
\nc{\bitbmc}{\begin{itemize}\begin{multicols}}
\nc{\bmc}{\begin{itemize}\begin{multicols}}
\nc{\emc}{\end{multicols}\end{itemize}\vspace{-3 truemm}}
\nc{\eit}{\end{itemize}\vspace{-3 truemm}}
\nc{\ben}{\vspace{-3 truemm}\begin{enumerate}}
\nc{\een}{\end{enumerate}\vspace{-3 truemm}}
\nc{\eitres}{\end{itemize}}

\nc{\set}[2]{\{ {#1} : {#2} \}} 
\nc{\bigset}[2]{\big\{ {#1}: {#2} \big\}} 
\nc{\Bigset}[2]{\left\{ \,{#1} :{#2}\, \right\}}

\nc{\pres}[2]{\la #1 : #2 \ra}
\nc{\bigpres}[2]{\big\la {#1} : {#2} \big\ra}
\nc{\Bigpres}[2]{\Big\la \,{#1}:{#2}\, \Big\ra}
\nc{\Biggpres}[2]{\Bigg\la {#1} : {#2} \Bigg\ra}

\nc{\pf}{\noindent{\bf Proof.}  }
\nc{\epf}{\hfill$\Box$\bigskip}
\nc{\epfres}{\hfill$\Box$}
\nc{\pfnb}{\pf}
\nc{\epfnb}{\bigskip}
\nc{\pfthm}[1]{\bigskip \noindent{\bf Proof of Theorem \ref{#1}}\,\,  } 
\nc{\pfprop}[1]{\bigskip \noindent{\bf Proof of Proposition \ref{#1}}\,\,  } 
\nc{\epfreseq}{\tag*{$\Box$}}

\nc{\uvert}[1]{\fill (#1,2)circle(.2);}
\rnc{\lvert}[1]{\fill (#1,0)circle(.2);}
\nc{\guvert}[1]{\fill[lightgray] (#1,2)circle(.2);}
\nc{\glvert}[1]{\fill[lightgray] (#1,0)circle(.2);}
\nc{\uvertx}[2]{\fill (#1,#2)circle(.2);}
\nc{\guvertx}[2]{\fill[lightgray] (#1,#2)circle(.2);}
\nc{\uvertxs}[2]{
\foreach \x in {#1}
{ \uvertx{\x}{#2}  }
}
\nc{\guvertxs}[2]{
\foreach \x in {#1}
{ \guvertx{\x}{#2}  }
}

\nc{\uvertth}[2]{\fill (#1,2)circle(#2);}
\nc{\lvertth}[2]{\fill (#1,0)circle(#2);}
\nc{\uvertths}[2]{
\foreach \x in {#1}
{ \uvertth{\x}{#2}  }
}
\nc{\lvertths}[2]{
\foreach \x in {#1}
{ \lvertth{\x}{#2}  }
}

\nc{\vertlabel}[2]{\draw(#1,2+.3)node{{\tiny $#2$}};}
\nc{\vertlabelh}[2]{\draw(#1,2+.4)node{{\tiny $#2$}};}
\nc{\vertlabelhh}[2]{\draw(#1,2+.6)node{{\tiny $#2$}};}
\nc{\vertlabelhhh}[2]{\draw(#1,2+.64)node{{\tiny $#2$}};}
\nc{\vertlabelup}[2]{\draw(#1,4+.6)node{{\tiny $#2$}};}
\nc{\vertlabels}[1]{
{\foreach \x/\y in {#1}
{ \vertlabel{\x}{\y} }
}
}

\nc{\dvertlabel}[2]{\draw(#1,-.4)node{{\tiny $#2$}};}
\nc{\dvertlabels}[1]{
{\foreach \x/\y in {#1}
{ \dvertlabel{\x}{\y} }
}
}
\nc{\vertlabelsh}[1]{
{\foreach \x/\y in {#1}
{ \vertlabelh{\x}{\y} }
}
}
\nc{\vertlabelshh}[1]{
{\foreach \x/\y in {#1}
{ \vertlabelhh{\x}{\y} }
}
}
\nc{\vertlabelshhh}[1]{
{\foreach \x/\y in {#1}
{ \vertlabelhhh{\x}{\y} }
}
}

\nc{\vertlabelx}[3]{\draw(#1,2+#3+.6)node{{\tiny $#2$}};}
\nc{\vertlabelxs}[2]{
{\foreach \x/\y in {#1}
{ \vertlabelx{\x}{\y}{#2} }
}
}

\nc{\vertlabelupdash}[2]{\draw(#1,2.7)node{{\tiny $\phantom{'}#2'$}};}
\nc{\vertlabelupdashess}[1]{
{\foreach \x/\y in {#1}
{\vertlabelupdash{\x}{\y}}
}
}

\nc{\vertlabeldn}[2]{\draw(#1,0-.6)node{{\tiny $\phantom{'}#2'$}};}
\nc{\vertlabeldnph}[2]{\draw(#1,0-.6)node{{\tiny $\phantom{'#2'}$}};}

\nc{\vertlabelups}[1]{
{\foreach \x in {#1}
{\vertlabel{\x}{\x}}
}
}
\nc{\vertlabeldns}[1]{
{\foreach \x in {#1}
{\vertlabeldn{\x}{\x}}
}
}
\nc{\vertlabeldnsph}[1]{
{\foreach \x in {#1}
{\vertlabeldnph{\x}{\x}}
}
}

\nc{\dotsup}[2]{\draw [dotted] (#1+.6,2)--(#2-.6,2);}
\nc{\dotsupx}[3]{\draw [dotted] (#1+.6,#3)--(#2-.6,#3);}
\nc{\dotsdn}[2]{\draw [dotted] (#1+.6,0)--(#2-.6,0);}
\nc{\dotsups}[1]{\foreach \x/\y in {#1}
{ \dotsup{\x}{\y} }
}
\nc{\dotsupxs}[2]{\foreach \x/\y in {#1}
{ \dotsupx{\x}{\y}{#2} }
}
\nc{\dotsdns}[1]{\foreach \x/\y in {#1}
{ \dotsdn{\x}{\y} }
}

\nc{\nodropcustpartn}[3]{
\begin{tikzpicture}[scale=.3]
\foreach \x in {#1}
{ \uvert{\x}  }
\foreach \x in {#2}
{ \lvert{\x}  }
#3 \end{tikzpicture}
}

\nc{\custpartn}[3]{{\lower1.4 ex\hbox{
\begin{tikzpicture}[scale=.3]
\foreach \x in {#1}
{ \uvert{\x}  }
\foreach \x in {#2}
{ \lvert{\x}  }
#3 \end{tikzpicture}
}}}

\nc{\smcustpartn}[3]{{\lower0.7 ex\hbox{
\begin{tikzpicture}[scale=.2]
\foreach \x in {#1}
{ \uvert{\x}  }
\foreach \x in {#2}
{ \lvert{\x}  }
#3 \end{tikzpicture}
}}}

\nc{\dropcustpartn}[3]{{\lower5.2 ex\hbox{
\begin{tikzpicture}[scale=.3]
\foreach \x in {#1}
{ \uvert{\x}  }
\foreach \x in {#2}
{ \lvert{\x}  }
#3 \end{tikzpicture}
}}}

\nc{\dropcustpartnx}[4]{{\lower#4 ex\hbox{
\begin{tikzpicture}[scale=.4]
\foreach \x in {#1}
{ \uvert{\x}  }
\foreach \x in {#2}
{ \lvert{\x}  }
#3 \end{tikzpicture}
}}}

\nc{\dropcustpartnxy}[3]{\dropcustpartnx{#1}{#2}{#3}{4.6}}

\nc{\uvertsm}[1]{\fill (#1,2)circle(.15);}
\nc{\lvertsm}[1]{\fill (#1,0)circle(.15);}
\nc{\vertsm}[2]{\fill (#1,#2)circle(.15);}

\nc{\bigdropcustpartn}[3]{{\lower6.93 ex\hbox{
\begin{tikzpicture}[scale=.6]
\foreach \x in {#1}
{ \uvertsm{\x}  }
\foreach \x in {#2}
{ \lvertsm{\x}  }
#3 \end{tikzpicture}
}}}

\nc{\gcustpartn}[5]{{\lower1.4 ex\hbox{
\begin{tikzpicture}[scale=.3]
\foreach \x in {#1}
{ \uvert{\x}  }
\foreach \x in {#2}
{ \guvert{\x}  }
\foreach \x in {#3}
{ \lvert{\x}  }
\foreach \x in {#4}
{ \glvert{\x}  }
#5 \end{tikzpicture}
}}}

\nc{\gcustpartndash}[5]{{\lower6.97 ex\hbox{
\begin{tikzpicture}[scale=.3]
\foreach \x in {#1}
{ \uvert{\x}  }
\foreach \x in {#2}
{ \guvert{\x}  }
\foreach \x in {#3}
{ \lvert{\x}  }
\foreach \x in {#4}
{ \glvert{\x}  }
#5 \end{tikzpicture}
}}}

\nc{\stline}[2]{\draw(#1,2)--(#2,0);}
\nc{\stlines}[1]{
{\foreach \x/\y in {#1}
{ \stline{\x}{\y} }
}
}

\nc{\uarcs}[1]{
{\foreach \x/\y in {#1}
{ \uarc{\x}{\y} }
}
}

\nc{\darcs}[1]{
{\foreach \x/\y in {#1}
{ \darc{\x}{\y} }
}
}

\nc{\stlinests}[1]{
{\foreach \x/\y in {#1}
{ \stlinest{\x}{\y} }
}
}

\nc{\stlineds}[1]{
{\foreach \x/\y in {#1}
{ \stlined{\x}{\y} }
}
}

\nc{\gstline}[2]{\draw[lightgray](#1,2)--(#2,0);}
\nc{\gstlines}[1]{
{\foreach \x/\y in {#1}
{ \gstline{\x}{\y} }
}
}

\nc{\gstlinex}[3]{\draw[lightgray](#1,2+#3)--(#2,0+#3);}
\nc{\gstlinexs}[2]{
{\foreach \x/\y in {#1}
{ \gstlinex{\x}{\y}{#2} }
}
}

\nc{\stlinex}[3]{\draw(#1,2+#3)--(#2,0+#3);}
\nc{\stlinexs}[2]{
{\foreach \x/\y in {#1}
{ \stlinex{\x}{\y}{#2} }
}
}

\nc{\darcx}[3]{\draw(#1,0)arc(180:90:#3) (#1+#3,#3)--(#2-#3,#3) (#2-#3,#3) arc(90:0:#3);}
\nc{\darc}[2]{\darcx{#1}{#2}{.4}}
\nc{\uarcx}[3]{\draw(#1,2)arc(180:270:#3) (#1+#3,2-#3)--(#2-#3,2-#3) (#2-#3,2-#3) arc(270:360:#3);}
\nc{\uarc}[2]{\uarcx{#1}{#2}{.4}}

\nc{\darcxx}[4]{\draw(#1,0+#4)arc(180:90:#3) (#1+#3,#3+#4)--(#2-#3,#3+#4) (#2-#3,#3+#4) arc(90:0:#3);}
\nc{\uarcxx}[4]{\draw(#1,2+#4)arc(180:270:#3) (#1+#3,2-#3+#4)--(#2-#3,2-#3+#4) (#2-#3,2-#3+#4) arc(270:360:#3);}

\makeatletter
\newcommand\footnoteref[1]{\protected@xdef\@thefnmark{\ref{#1}}\@footnotemark}
\makeatother

\newcounter{theorem}
\numberwithin{theorem}{section}

\newtheorem{thm}[theorem]{Theorem}
\newtheorem{lemma}[theorem]{Lemma}
\newtheorem{cor}[theorem]{Corollary}
\newtheorem{prop}[theorem]{Proposition}

\theoremstyle{definition}

\newtheorem{rem}[theorem]{Remark}
\newtheorem{defn}[theorem]{Definition}
\newtheorem{eg}[theorem]{Example}
\newtheorem{ass}[theorem]{Assumption}

\title{Presentations for rook partition monoids and algebras \\ and their singular ideals}

\date{}

\author{
James East\\
{\footnotesize \emph{Centre for Research in Mathematics; School of Computing, Engineering and Mathematics}}\\
{\footnotesize \emph{University of Western Sydney, Locked Bag 1797, Penrith NSW 2751, Australia}}\\
{\footnotesize {\tt J.East\,@\,uws.edu.au}}
}

\maketitle

\vspace{-0.5cm}

\begin{abstract}
We obtain several presentations by generators and relations for the rook partition monoids and algebras, as well as their singular ideals.  Among other results, we also calculate the minimal sizes of generating sets (some of our presentations use such minimal-size generating sets), and show that the singular part of the rook partition monoid is generated by its idempotents.

{\it Keywords}: Partition algebra, Partition monoid, Rook partition algebra, Rook partition monoid, Singular ideal, Presentations, Rank, Idempotent rank.

MSC: 20M05; 20M20.
\end{abstract}

\section{Introduction}\label{sect:intro}

\emph{Diagram algebras} have diverse origins and applications.  For example, see Brauer \cite{Brauer1937} on invariant theory; Jones \cite{Jones1987} and Kauffman \cite{Kauffman1990} on knot theory; Temperley and Lieb \cite{TL1971}, Jones \cite{Jones1994_2} and Martin~\cite{Martin1994} on statistical mechanics; and more.  One unifying theme is that diagram algebras often occur as \emph{centraliser algebras} of classical groups, leading to various interesting extensions of (classical) \emph{Schur-Weyl duality} \cite{Weyl_book}.  For example, the \emph{Brauer algebras} are related in this way to the orthogonal groups \cite{Brauer1937}, and the \emph{partition algebras} to symmetric groups \cite{Martin1994,Jones1994_2}.
In a highly influential~2005 paper, Halverson and Ram \cite{HR2005} gave (among many other things) an account of Schur-Weyl duality in the case of the partition algebras, making crucial use of a \emph{tower of algebras} that saw the partition algebras $\C A_k(n)$ embedded into $\C A_{k+1}(n)$ via an intermediate subalgebra: 
\[
\cdots\hookrightarrow\C A_k(n)\hookrightarrow\C A_{k+\frac12}(n)\sub\C A_{k+1}(n)\hookrightarrow\cdots.
\]
Halverson and Ram \cite{HR2005} attributed their understanding of the ``existence and importance'' of the algebras $\C A_{k+\frac12}(n)$ to Cheryl Grood, who studied them in their own right in \cite{Grood06}, where they were called \emph{rook partition algebras}, and given their own diagrammatic interpretation (see Section \ref{subsect:RPn} below for details); Grood also noted that these intermediate algebras were used in earlier work of Martin\mbox{\cite{Martin1996,Martin2000}}.  The reason for the name is due to a connection with the so-called \emph{rook monoids} (and associated algebras and deformations) studied by Halverson, Solomon and others \cite{DHP2003,Halverson2004,Grood2002,HR2001,Paget2006,Solomon2002,JEcais}.  As noted by Grood \cite{Grood06}, Solomon's discovery \cite{Solomon2002} of a Schur-Weyl duality for rook monoid algebras (see also \cite{KM2008}) led to the investigation of a number of other ``rook diagram algebras'', such the \emph{rook Brauer algebras} \cite{MarMaz2014,HD2014}, \emph{Motzkin algebras} \cite{BH2014} and more.  Such studies, and other considerations often to do with representation theory and/or statistical mechanics, have led to the discovery and investigation of a great many other families of
diagram algebras \cite{CGM2003,DO2014,Bloss2003,PK2004,Orellana2007,Kennedy2007,MW2000,MS1994}.

In this article, we continue the study of the rook partition algebras, with our goal being to derive \emph{presentations} by generators and relations.  This continues a theme initiated by the current author in \cite{JEgrpm}, where such presentations were obtained for the partition algebras themselves.  Presentations are extremely useful tools for algebraists: for one thing, they allow representations (homomorphisms into other algebraic structures) to be defined by specifying the images of the generators and checking that the relations are preserved; this is especially helpful when the algebra in question is as complicated as the partition algebra.
With respect to the partition algebras in particular, we refer to the recent work of Enyang on Jucys-Murphy elements \cite{Enyang2013_1} and seminormal forms \cite{Enyang2013_2}, in which the presentations from \cite{JEgrpm,HR2005} played a crucial role.  Presentations also feature heavily in the work of Lehrer and Zhang on diagram categories and invariant theory \cite{LZ2015,LZ2012}.

Key to the approach used in \cite{JEgrpm} was the observation of Wilcox \cite{Wilcox2007} (also implicit in \cite{HR2005}) that diagram algebras (including the partition algebras) arise as \emph{twisted semigroup algebras} of corresponding \emph{diagram semigroups}.  This allows one to obtain information (concerning cellularity \cite{Wilcox2007,EastGray,DEG2015,GL1996} or presentations \cite{JEpnsn,JEgrpm}, for example) about the algebras from corresponding information about the associated semigroups.  Conversely, the theory of diagram algebras has led to a number of important families of semigroups and monoids that have been studied with increasing vigour in recent years; see for example \cite{PHY2013,FL2011,LF2006,BDP2002,JEgrpm,JEpnsn,DE2015,DE2016,DEEFHHL1,DEEFHHL2,DEG2015,EF,EastGray,Maz2002,Maz1998, KM2006,Maltcev2007}, and especially the work of Auinger and his collaborators \cite{ACHLV2015,ADV2012_2,ADV2012,Auinger2014} on equational theories of involution semigroups.  A number of other authors have studied presentations of diagram semigroups \cite{BDP2002,PHY2013,Maltcev2007,KM2006}.  We hope that the techniques we introduce here to deal with the more complicated \emph{rook} partition algebras will prove useful in other investigations: for example, in the context of \emph{quasi-partition algebras} \cite{DO2014} or \emph{coloured partition algebras} \cite{Bloss2003,PK2004} and their rook versions \cite{KM2013}.

As stated above, our focus in this article is on the rook partition monoid and algebra, $\RP_n$ and $F^\tau[\RP_n]$, and their singular ideals, $\RPnSn$ and $F^\tau[\RPnSn]$; see Sections \ref{sect:prelim} and~\ref{sect:algebras} for precise definitions.  In particular, we obtain a number of presentations for each of these algebraic systems (Theorems \ref{thm:RPnSn}, \ref{thm:RPn}, \ref{thm:RPn2}, \ref{thm:RPn3} and \ref{thm:algebras}); among other results, we calculate the smallest sizes of generating sets (Theorems \ref{thm:rankRPnSn} and \ref{thm:rankRPn}), and show that $\RPnSn$ is generated by its idempotents (Proposition~\ref{prop:RPnSn_gen} and Remark \ref{rem:RPnSn_gen}).  Our approach is quite different to previous studies of (singular) diagram semigroups and algebras \cite{JEgrpm,JEpnsn,KM2006,Maltcev2007,BDP2002,PHY2013}, in the sense that we first obtain a presentation for the singular rook partition monoid $\RPnSn$ (the subject of Section \ref{sect:RPnSn}), and then use this to bootstrap up to a number of presentations for the (full) rook partition monoid $\RP_n$ (in Section~\ref{sect:RPn}).  
(Although this appears to be the first instance of such a ``singular first'' approach, we note that a somewhat similar method was used in the author's recent work on (singular) symmetric inverse semigroups in \cite{JEinsn2}.)  
Our approach makes crucial use of the author's presentations for the (ordinary) partition monoid $\P_n$ \cite{JEgrpm} and its singular ideal $\PnSn$ \cite{JEpnsn}; these are stated in Section \ref{sect:prelim}, along with various definitions, notations, background information, and illustrative examples.  Finally, in Section \ref{sect:algebras}, we apply general results on twisted semigroup algebras from \cite{JEgrpm} to obtain (algebra) presentations for the rook partition algebra $F^\tau[\RP_n]$ and its singular ideal $F^\tau[\RPnSn]$; see Theorem \ref{thm:algebras}.

\section{Preliminaries}\label{sect:prelim}

In this section, we gather the preliminary material we will need in our investigations.  In particular, we define the rook partition monoids, identify several key submonoids, develop the notation and parameters needed to formulate and prove our results, and review previous results on (ordinary) partition monoids that will play a role in our study.  (We postpone the definition of the rook partition \emph{algebra} until Section \ref{sect:algebras}.)

\subsection{The rook partition monoid}\label{subsect:RPn}

We adapt the treatment of Grood \cite{Grood06}.  Fix a non-negative integer $n$, and write ${\bn=\{1,\ldots,n\}}$ and $\bn'=\{1',\ldots,n'\}$.  By a \emph{rook partition of degree $n$}, we mean a set partition of a subset of $\bn\cup\bn'$: i.e., a collection $\al=\set{A_i}{i\in I}$, for some (possibly empty) indexing set $I$, where the $A_i$ are pairwise disjoint non-empty subsets of $\bn\cup\bn'$; the $A_i$ are called the \emph{blocks} of~$\al$.  
We write $\supp(\al)=\bigcup_{i\in I}A_i$, and call this the \emph{support} of $\al$.  The elements of $(\bn\cup\bn')\sm\supp(\al)$ are called the \emph{rook dots} of $\al$.  We write $\RP_n$ for the set of all rook partitions of degree $n$.  For example, the rook partition
\[
\al=\big\{ \{1,2,4,3'\}, \{5,6,4',5'\}, \{7,8,8'\}, \{2',6',7'\}, \{9',10'\} \big\}\in\RP_{10}
\]
has support $\{1,2,4,5,6,7,8,2',3',4',5',6',7',8',9',10'\}$, and rook dots $3,9,10,1'$.

\begin{figure}[H]
\begin{center}
\begin{tikzpicture}[scale=.39]
\begin{scope}[shift={(0,0)}]	
\uvws{3,9,10}
\lvws{1}
\uverts{1,2,4,5,6,7,8}
\lverts{2,3,4,5,6,7,8,9,10}
\uarcs{1/2,2/4,5/6,7/8}
\darcs{4/5,6/7,9/10}
\darcx26{.8}
\stlines{4/3,5/5,8/8}
\draw(0.5,1)node[left]{$\al=$};
\end{scope}
\begin{scope}[shift={(0,-4)}]	
\uvws{1,8,10}
\lvws{3}
\uverts{2,...,7,9}
\lverts{1,2,4,5,...,10}
\uarcs{2/3,5/6}
\darcs{1/2,4/5,6/7,8/9}
\darcx7{10}{.8}
\stlines{4/4,7/7,9/9}
\draw(0.5,1)node[left]{$\be=$};
\draw(13,3)node{$\longrightarrow$};
\end{scope}
\begin{scope}[shift={(15,-1)}]	
\uverts{1,2,4,5,6,7,8}
\lverts{2,3,4,5,6,7,8,9,10}
\uarcs{1/2,2/4,5/6,7/8}
\darcs{4/5,6/7,9/10}
\darcx26{.8}
\stlines{4/3,5/5,8/8}
\uvws{3,9,10}
\lvws{1}
\end{scope}
\begin{scope}[shift={(15,-3)}]	
\uvws{1,8,10}
\lvws{3}
\uverts{2,...,7,9}
\lverts{1,2,4,5,...,10}
\uarcs{2/3,5/6}
\darcs{1/2,4/5,6/7,8/9}
\darcx7{10}{.8}
\stlines{4/4,7/7,9/9}
\draw(13,2)node{$\longrightarrow$};
\end{scope}
\begin{scope}[shift={(30,-2)}]	
\uvws{3,7,8,9,10}
\lvws{3,8,9}
\uverts{1,2,4,5,6}
\lverts{1,2,4,5,6,7,10}
\uarcs{1/2,2/4,4/5,5/6}
\darcs{1/2,4/5,5/6,6/7,7/10}
\stlines{5/5}
\draw(10.5,1)node[right]{$=\al\be$};
\end{scope}
\end{tikzpicture}
\end{center}
\vspace{-5mm}
\caption{Two rook partitions $\al,\be\in\RP_{10}$ (left), the product $\al\be\in\RP_{10}$ (right), and the product graph $\Ga(\al,\be)$ (centre).}
\label{fig:RP10}
\end{figure}
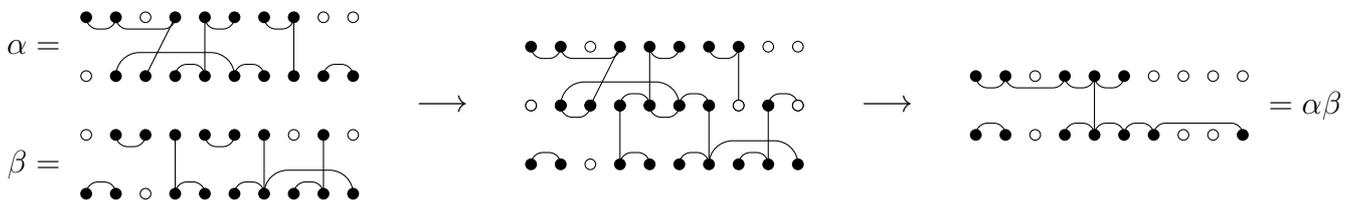

We may represent an element $\al\in\RP_n$ graphically as follows.  We arrange vertices $1,\ldots,n$ in a horizontal row (increasing from left to right) with vertices $1',\ldots,n'$ directly below.  We colour each vertex from $\supp(\al)$ black, and each rook dot white, and add edges between black vertices in such a way that two vertices 
are joined by a path if and only if they belong to the same block of $\al$.  Such a graph is called a \emph{rook diagram} of $\al$.  The rook partition $\al\in\RP_{10}$ defined at the end of the previous paragraph is pictured in Figure \ref{fig:RP10} (top left).
The graphical representation of a rook partition is not unique, but, as in \cite{Grood06}, we regard two rook diagrams as \emph{equivalent} if they have the same rook dots and the same connected components.  Generally, we identify a rook partition $\al\in\RP_n$ with any rook diagram representing it.

To describe the product in $\RP_n$, consider two rook partitions $\al,\be\in\RP_n$.  Write $\bn''=\{1'',\ldots,n''\}$.  Let $\alb$ be the graph obtained from $\al$ by changing the label of each lower vertex $i'$ to~$i''$.  Similarly, let $\beb$ be the graph obtained from $\be$ by changing the label of each upper vertex~$i$ to~$i''$.  Consider now the graph $\Ga(\al,\be)$ on the vertex set~$\bn\cup \bn'\cup \bn''$ obtained by joining $\alb$ and~$\beb$ together so that each lower vertex $i''$ of $\alb$ is identified with the corresponding upper vertex $i''$ of $\beb$, and colouring a vertex $i''$ white 
in $\Ga(\al,\be)$ if it is white in either $\alb$ or $\beb$.
We call $\Ga(\al,\be)$ the \emph{product graph of $\al,\be$}.  We now let $\al\be\in\RP_n$ be the (rook partition corresponding to any) rook diagram such that:
\bit
\item $\al\be$ has a rook dot $x\in\bn\cup\bn'$ if and only if $x$ is connected by a (possibly empty) path to a white vertex in $\Ga(\al,\be)$, and
\item two black vertices $x,y\in\bn\cup\bn'$ are connected by a path in $\al\be$ if and only if they are connected by a path in $\Ga(\al,\be)$.
\eit
An example calculation (with $n=10$) is carried out in Figure \ref{fig:RP10}.

This operation is associative, and gives $\RP_n$ the structure of a monoid; the identity element is
\[
1 = \big\{\{1,1'\},\ldots,\{n,n'\}\big\} = \custpartn{1,2,4}{1,2,4}{\dotsups{2/4}\dotsdns{2/4}\stlines{1/1,2/2,4/4}\vertlabelshh{1/1,4/n}} \in\RP_n.
\]
Before we can say more about the structure of $\RP_n$, we first introduce some notation.  Let $\al\in\RP_n$.  We call a block of $\al$ a \emph{transversal block} if it has non-empty intersection with both $\bn$ and $\bn'$, and a \emph{non-transversal block} otherwise; we will also distinguish between \emph{upper} and \emph{lower} non-transversal blocks (defined in an obvious way).  The \emph{rank} of $\al$, denoted $\rank(\al)$, is defined to be the number of transversal blocks of $\al$.  Note that $0\leq\rank(\al)\leq n$.  For $x\in\supp(\al)$, we write $[x]_\al$ for the block of $\al$ containing $x$; we also define $[x]_\al=\{x\}$ if $x$ is a rook dot of $\al$.
We then define the \emph{domain} and \emph{codomain} of $\al$ to be the sets
\begin{align*}
\dom(\al) = \set{x\in\bn}{[x]_\al\cap\bn'\not=\emptyset} &\AND
\codom(\al) = \set{x\in\bn}{[x']_\al\cap\bn\not=\emptyset}.
\intertext{We also define the \emph{kernel} and \emph{cokernel} of $\al$ to be the equivalences}
\ker(\al) = \bigset{(x,y)\in\bn\times\bn}{[x]_\al=[y]_\al} &\AND
\coker(\al) = \bigset{(x,y)\in\bn\times\bn}{[x']_\al=[y']_\al}.
\end{align*}
For example, with $\al\in\RP_{10}$ as in Figure \ref{fig:RP10} (top left), 
\begin{gather*}
\rank(\al)=3\COMMA \dom(\al)=\{1,2,4,5,6,7,8\} \COMMA \codom(\al) = \{3,4,5,8\}, \\
\ker(\al) = (\ 1,2,4\ | \ 3 \ |\ 5,6\ |\ 7,8 \ |\ 9 \ |\ 10\ ) \COMMA \coker(\al) = (\ 1 \ |\ 2,6,7 \ |\ 3 \ |\ 4,5 \ |\ 8 \ |\ 9,10\ ),
\end{gather*}
using an obvious notation for equivalences.

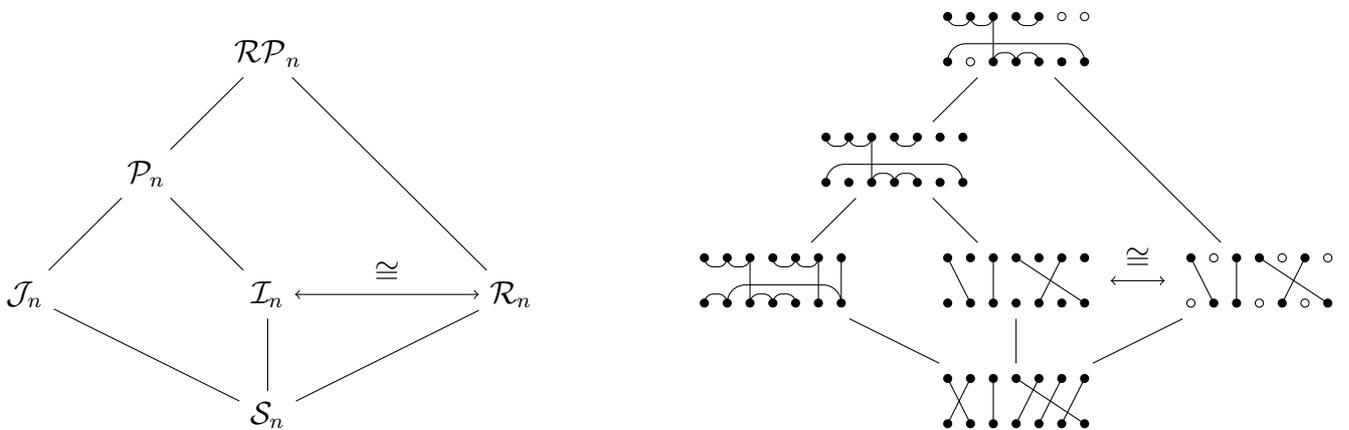
\begin{figure}[h]
\begin{center}
\begin{tikzpicture}[scale=.8]
\node (A) at (4,6) {$\RP_n$};
\node (B) at (2,4) {$\P_n$};
\node (C) at (0,2) {$\J_n$};
\node (D) at (4,2) {$\I_n$};
\node (E) at (8,2) {$\R_n$};
\node (F) at (4,0) {$\S_n$};
\draw (F)--(C)--(B)--(D)--(F)--(E)--(A)--(B);
\draw[<->] (D)-- node[above]{$\cong$} (E);
\end{tikzpicture}
\qquad\qquad
\begin{tikzpicture}[scale=.8]
\node (A) at (4,6) {$\custpartn{1,2,3,4,5}{1,3,4,5,6,7}{\uvw6\uvw7\lvw2\stlines{3/3}\uarcs{1/2,2/3,4/5}\darcs{3/4,4/5}\darcx17{.8}}$};
\node (B) at (2,4) {$\custpartn{1,...,7}{1,...,7}{\stlines{3/3}\uarcs{1/2,2/3,4/5}\darcs{3/4,4/5}\darcx17{.8}}$};
\node (C) at (0,2) {$\custpartn{1,...,7}{1,...,7}{\stlines{3/3,6/6,7/7}\uarcs{1/2,2/3,4/5,5/6}\darcs{1/2,3/4,4/5}\darcx27{.8}}$};
\node (D) at (4,2) {$\custpartn{1,...,7}{1,...,7}{\stlines{1/2,3/3,4/7,6/5}}$};
\node (E) at (8,2) {$\custpartn{1,3,4,6}{2,3,5,7}{\uvws{2,5,7}\lvws{1,4,6}\stlines{1/2,3/3,4/7,6/5}}$};
\node (F) at (4,0) {$\custpartn{1,...,7}{1,...,7}{\stlines{1/2,2/1,3/3,4/7,7/6,6/5,5/4}}$};
\draw (F)--(C)--(B)--(D)--(F)--(E)--(A)--(B);
\draw[<->] (D)-- node[above]{$\cong$} (E);
\end{tikzpicture}
\end{center}
\vspace{-5mm}
\caption{Important submonoids of $\RP_n$ (left) and representative elements from each submonoid (right).}
\label{fig:submonoids}
\end{figure}

It is immediate from the definitions that the following hold for all $\al,\be\in\RP_n$:
\[
\begin{array}{rclcrcl}
\dom(\al\be) \sub \dom(\al) \COMMA
\ker(\al\be)\sp \ker(\al) \COMMA
\codom(\al\be) \sub \codom(\be) \COMMA
\coker(\al\be)\sp \coker(\be).
\end{array}
\]
We write $\De=\set{(x,x)}{x\in\bn}$ for the trivial relation on $\bn$: i.e., the equality relation.
The rook partition monoid $\RP_n$ contains a number of important submonoids, 
defined as follows:
\bit
\item $\P_n=\set{\al\in\RP_n}{\supp(\al)=\bn\cup\bn'}$, the \emph{partition monoid} \cite{HR2005,JEgrpm};
\item $\I_n=\set{\al\in\P_n}{\ker(\al)=\coker(\al)=\De}$, the \emph{symmetric inverse monoid} \cite{Lipscombe1996,Lawson1998};
\item $\J_n=\set{\al\in\P_n}{\dom(\al)=\codom(\al)=\bn}$, the \emph{dual symmetric inverse monoid} \cite{FL1998};
\item $\S_n = \set{\al\in\RP_n}{\rank(\al)=n}$, the \emph{symmetric group} \cite{Cameron1999}; and
\item $\R_n = \set{\al\in\RP_n}{\ker(\al)=\coker(\al)=\De,\ \supp(\al)=\dom(\al)\cup\codom(\al)'}$, the \emph{rook monoid}~\cite{Grood2002,Solomon2002}.
\eit
Elements of these sumbonoids are pictured in Figure \ref{fig:submonoids}, along with the various containments among the submonoids.  

Note that $\I_n$ and $\R_n$ are isomorphic (via an obvious map $\I_n\to\R_n$ that colours all non-transversal blocks of $\al\in\I_n$ white), and that $\S_n$ is the group of units of $\RP_n$ (in fact, the group of units of all the stated submonoids).  It follows from 
\cite[Theorem 32]{JEgrpm} that $\P_n$ is equal to the \emph{join} $\I_n\vee\J_n=\la\I_n\cup\J_n\ra$.  The join $\J_n\vee\R_n$ is the \emph{partial dual symmetric inverse monoid} studied in \cite{KMU2015,KM2011}.  The join $\I_n\vee\R_n$ is a ``rook version'' of the symmetric inverse monoid; to the author's knowledge, it has not been explicitly studied in the literature.  As noted in \cite{Grood06}, $\RP_n$ is isomorphic to the submonoid of $\P_{n+1}$ consisting of all (ordinary) partitions of degree $n+1$ such that $n+1$ and $(n+1)'$ belong to the same block; these submonoids were introduced in \cite{Martin1996,Martin2000} and played a central role in \cite{HR2005}.

We now describe a convenient notation for rook partitions, extending the notation introduced for ordinary partitions in \cite{EF}.  With this in mind, let $\al\in\RP_n$.  We write
\[
\al = \left( \begin{array}{c|c|c|c|c|c||c} 
\!\! A_1 \ & \cdots & A_r & C_1 & \cdots & C_p & P \ \ \\ \cline{4-7}
\!\! B_1 \ & \cdots & B_r & D_1 & \cdots & D_q & Q \ \
\end{array} \!\!\! \right)
\]
to indicate that $\al$ has (writing $A'=\set{a'}{a\in A}$ for $A\sub\bn$):
\bit
\item transversal blocks $A_i\cup B_i'$, for each $1\leq i\leq r$, 
\item upper non-transversal blocks $C_i$, for each $1\leq i\leq p$, 
\item lower non-transversal blocks $D_i'$, for each $1\leq i\leq q$, and 
\item rook dots $x$, for each $x\in P\cup Q'$.
\eit
So, for example, with $\al\in\RP_{10}$ as in Figure \ref{fig:RP10} (top left), we have
\[
\al= \left( \begin{array}{c|c|c|c|c||c} 
\!\! 1,2,4 & 5,6 & 7,8 & \multicolumn{2}{c||}{\ }  & 3,9,10 \ \ \\ \cline{4-6}
\!\! 3 & 4,5 & 8 & 2,6,7 & 9,10 & 1 \ \
\end{array} \!\!\! \right).
\]
In the above notation, it is possible that any of $r,p,q,|P|,|Q|$ could be $0$, in which case we may use simplified versions of the above notation.  If $\al$ has no rook dots (so that $\al\in\P_n$), we will omit $P$ and $Q$, and write
\[
\al = \partnlist{A_1}{A_r}{C_1}{C_p}{B_1}{B_r}{D_1}{D_q},
\]
as in \cite{EF}.
Similarly, if $\al$ has no rook dots and no non-transversal blocks (so that $\al\in\J_n$), then we will simply write
\[
\al=\partnlistnodefshort{A_1}{A_r}{B_1}{B_r}.
\]
If $\al\in\I_n$ (so $\ker(\al)=\coker(\al)=\De$ and $\al$ has no rook dots), we will write 
\[
\al=\partnlistnodefshort{a_1}{a_r}{b_1}{b_r} 
\]
to indicate that the transversal blocks of $\al$ are $\{a_i,b_i'\}$ (for each $1\leq i\leq r$) and that every other block is a (non-rook) singleton.  We will use other simplifications from time to time, but it will always be clear what is meant.

Finally, we note that $\RP_n$ has a natural (anti-)involution ${}^*:\RP_n\to\RP_n:\al\mt\al^*$, defined by
\[
\left( \begin{array}{c|c|c|c|c|c||c} 
\!\! A_1 \ & \cdots & A_r & C_1 & \cdots & C_p & P \ \ \\ \cline{4-7}
\!\! B_1 \ & \cdots & B_r & D_1 & \cdots & D_q & Q \ \
\end{array} \!\!\! \right)^*
=
\left( \begin{array}{c|c|c|c|c|c||c} 
\!\! B_1 \ & \cdots & B_r & D_1 & \cdots & D_q & Q \ \ \\ \cline{4-7}
\!\! A_1 \ & \cdots & A_r & C_1 & \cdots & C_p & P \ \
\end{array} \!\!\! \right).
\]
Diagrammatically, $\al^*$ is obtained by turning (a graphical representation of) $\al$ upside-down.  This involution reflects the structure of $\RP_n$ as a \emph{regular $*$-semigroup} (as defined in \cite{NS1978}).  That is, the following hold for all $\al,\be\in\RP_n$ (as may easily be checked diagrammatically, or follows from the above-mentioned embedding $\RP_n\to\P_{n+1}$):
\[
(\al^*)^*=\al \COMMA (\al\be)^*=\be^*\al^* \COMMA \al\al^*\al=\al \COMMA \al^*\al\al^*=\al^*.
\]
This $*$-regular structure leads to a duality that will help simplify several proofs, and has played a very important role in many other studies of diagram monoids; see \cite{ACHLV2015,ADV2012_2,Auinger2014,EF,JEgrpm,JEpnsn,EastGray}, among many others.

We conclude this section with two structural results, concerning normal forms for $\P_n$ (Proposition~\ref{prop:Pn_normalform}) and $\RP_n$ (Proposition \ref{prop:RPn_normalform}).  They will both be useful on a number of occasions.  
The first of these gives a convenient factorisation for the elements of $\P_n$ that was used in the proof of \cite[Theorem~30]{EF}; its proof is easy and is omitted.

\ms
\begin{prop}\label{prop:Pn_normalform}
Let $\al\in\P_n$, and write
\[
\al = \partnlist{A_1}{A_r}{C_1}{C_p}{B_1}{B_r}{D_1}{D_q}.
\]
For each $i\in\br$, choose some $a_i\in A_i$ and $b_i\in B_i$.
Then $\al=\be\ga\de$, where
\[\epfreseq
\be = \partnlistnodef{A_1}{A_r}{C_1}{C_p}{A_1}{A_r}{C_1}{C_p} 
\COMMA
\ga = \partnlistnodefshort{a_1}{a_r}{b_1}{b_r}
\COMMA
\de = \partnlistnodef{B_1}{B_r}{D_1}{D_q}{B_1}{B_r}{D_1}{D_q}
.
\]
\end{prop}

\ms
\begin{rem}
Since $\be,\de\in\J_n$ and $\ga\in\I_n$ (as defined above), Proposition~\ref{prop:Pn_normalform} gives another proof of the above-mentioned fact that $\P_n=\I_n\vee\J_n=\la\I_n\cup\J_n\ra$.
\end{rem}

As an example to illustrate Proposition \ref{prop:Pn_normalform}, consider the partition
\begin{align*}
\al &= \big\{ \{1,2,4,3'\}, \{5,6,4',5'\}, \{7,8,8'\}, \{3\}, \{9,10\}, \{1'\}, \{2',6',7'\}, \{9',10'\} \big\} \in\P_{10}.
\end{align*}
Its factorisation $\al=\be\ga\de$ (for some choice of the $a_i,b_i$), as in Proposition \ref{prop:Pn_normalform}, is shown in Figure \ref{fig:a=bcd}.

\begin{figure}[H]
\begin{center}
\begin{tikzpicture}[scale=.4]
\begin{scope}[shift={(0,0)}]	
\uverts{1,...,10}
\lverts{1,...,10}
\uarcs{1/2,2/4,5/6,7/8,9/10}
\darcs{4/5,6/7,9/10}
\darcx26{.8}
\stlines{4/3,5/5,8/8}
\draw(0.5,1)node[left]{$\al=$};
\end{scope}
\begin{scope}[shift={(14,2)}]	
\uverts{1,...,10}
\lverts{1,...,10}
\uarcs{1/2,2/4,5/6,7/8,9/10}
\darcs{1/2,2/4,5/6,7/8,9/10}
\stlines{4/4,3/3,5/5,8/8,9/9}
\draw[|-|] (11.2,0)--(11.2,2);
\draw(11.2,1)node[right]{$\be$};
\end{scope}
\begin{scope}[shift={(14,0)}]	
\uverts{1,...,10}
\lverts{1,...,10}
\stlines{4/3,5/5,8/8}
\draw(-1,1)node[left]{$=$};
\draw[|-] (11.2,0)--(11.2,2);
\draw(11.2,1)node[right]{$\ga$};
\end{scope}
\begin{scope}[shift={(14,-2)}]	
\uverts{1,...,10}
\lverts{1,...,10}
\uarcs{4/5,6/7,9/10}
\uarcx26{.8}
\darcs{4/5,6/7,9/10}
\darcx26{.8}
\stlines{1/1,2/2,3/3,5/5,8/8,9/9}
\draw[|-] (11.2,0)--(11.2,2);
\draw(11.2,1)node[right]{$\de$};
\end{scope}
\end{tikzpicture}
\end{center}
\vspace{-5mm}
\caption{An illustration of the factorisation $\al=\be\ga\de$ from Proposition \ref{prop:Pn_normalform}, where $\al\in\P_{10}$.}
\label{fig:a=bcd}
\end{figure}

Consider a subset $A\sub\bn$, and write $A^c=\{i_1,\ldots,i_r\}$.  Here are elsewhere, we will write $A^c$ for the complement $\bn\sm A$.  We define \[
\ob_A = \left( \begin{array}{c|c|c||c} 
\!\! i_1 & \cdots & i_r & A \ \ \\ \cline{4-4}
\!\! i_1 & \cdots & i_r & A \ \
\end{array} \!\!\! \right)
\]
to be the (unique) element of $\RP_n$ with rook dots $x$ for each $x\in A\cup A'$, and transversal blocks~$\{j,j'\}$ for each $j\in A^c$.  (The reason for the (over-line) notation will become clear shortly.)  So, for example, with $n=10$ and $A=\{2,4,7,8,10\}$,
\[
\ob_A = \custpartn{1,3,5,6,9}{1,3,5,6,9}{\uvws{2,4,7,8,10}\lvws{2,4,7,8,10}\stlines{1/1,3/3,5/5,6/6,9/9}} \in\RP_{10}.
\]
To simplify the statement of the next result, we will allow ourselves to refer to the empty set $\emptyset$ as a block of any element of $\P_n$.

\ms
\begin{prop}\label{prop:RPn_normalform}
For any rook partition $\al\in\RP_n$, $\al=\ob_P\be\ob_Q$ for unique (possibly empty) subsets $P,Q\sub\bn$ and a unique partition $\be\in\P_n$ such that $P$ and $Q'$ are (possibly empty) blocks of $\be$.  Moreover, if $\al\in\RPnSn$, then $\be\in\PnSn$.
\end{prop}

\pf Write 
\[
\al = \left( \begin{array}{c|c|c|c|c|c||c} 
\!\! A_1 & \cdots & A_r & C_1 & \cdots & C_p & P \ \ \\ \cline{4-7}
\!\! B_1 & \cdots & B_r & D_1 & \cdots & D_q & Q \ \
\end{array} \!\!\! \right).
\]
Let $\be\in\P_n$ be the partition obtained by replacing the rook dots from $\al$ with the (possibly empty) blocks $P$ and $Q'$.  It is then clear that $\al=\ob_P\be\ob_Q$.  This establishes the existence part of the result.  For the uniqueness part, suppose that $\al=\ob_A\ga\ob_B$ where $A$ and $B'$ are (possibly empty) blocks of $\ga$.  The rook dots of $\ob_A\ga\ob_B$ are clearly the elements of $A\cup B'$.  The blocks of $\ga$ other than $A$ and $B'$ also clearly coincide with the (non-rook) blocks of $\ob_A\ga\ob_B$.  This all shows that $A=P$, $B=Q$ and $\ga=\be$.

Finally, suppose $\al\in\RPnSn$.  If $\al$ has a rook dot, then $\be$ has a non-transversal block.  If $\al$ has no rook dots, then $\be=\al$.  So in either case, $\be\in\PnSn$. \epf

\ms
\begin{rem}
An alternative factorisation $\al=\ob_P\be\ob_Q$ could be proved, with $P\cup Q'$ a (possibly empty) block of~$\be$ (instead of separate blocks $P,Q'$).  However, the implication ${\al\in\RPnSn \implies \be\in\PnSn}$ would not hold in general: namely, if $\al\in\R_n$ and $\rank(\al)=n-1$, then $\be\in\S_n$.
\end{rem}

As an example to illustrate Proposition \ref{prop:RPn_normalform}, consider (once again) the rook partition
\begin{align*}
\al &= \big\{ \{1,2,4,3'\}, \{5,6,4',5'\}, \{7,8,8'\}, \{2',6',7'\}, \{9',10'\} \big\} \in\RP_{10}.
\end{align*}
Its factorisation $\al=\ob_P\be\ob_Q$, as in Proposition \ref{prop:RPn_normalform}, is shown in Figure \ref{fig:a=obo}.

\begin{figure}[h]
\begin{center}
\begin{tikzpicture}[scale=.4]
\begin{scope}[shift={(0,0)}]	
\uvws{3,9,10}
\lvws{1}
\uverts{1,2,4,5,6,7,8}
\lverts{2,3,4,5,6,7,8,9,10}
\uarcs{1/2,2/4,5/6,7/8}
\darcs{4/5,6/7,9/10}
\darcx26{.8}
\stlines{4/3,5/5,8/8}
\draw(0.5,1)node[left]{$\al=$};
\end{scope}
\begin{scope}[shift={(14,3.5)}]	
\foreach \x in {1,...,10} {\draw[dotted](\x,-.2)--(\x,-1.5);}
\stlines{1/1,2/2,4/4,5/5,6/6,7/7,8/8}
\uvws{3,9,10}
\lvws{3,9,10}
\uverts{1,2,4,5,6,7,8}
\lverts{1,2,4,5,6,7,8}
\draw[|-|] (11.2,0)--(11.2,2);
\draw(11.2,1)node[right]{$\ob_P$};
\end{scope}
\begin{scope}[shift={(14,0)}]	
\uarcs{1/2,2/4,5/6,7/8,9/10}
\darcs{4/5,6/7,9/10}
\uarcx39{.8}
\darcx26{.8}
\stlines{4/3,5/5,8/8}
\uverts{1,...,10}
\lverts{1,...,10}
\draw(-1,1)node[left]{$=$};
\draw[|-|] (11.2,0)--(11.2,2);
\draw(11.2,1)node[right]{$\be$};
\foreach \x in {1,...,10} {\draw[dotted](\x,0)--(\x,-1.3);}
\end{scope}
\begin{scope}[shift={(14,-3.5)}]	
\stlines{2/2,3/3,4/4,5/5,6/6,7/7,8/8,9/9,10/10}
\uvws{1}
\lvws{1}
\uverts{2,3,4,5,6,7,8,9,10}
\lverts{2,3,4,5,6,7,8,9,10}
\draw[|-|] (11.2,0)--(11.2,2);
\draw(11.2,1)node[right]{$\ob_Q$};
\end{scope}
\end{tikzpicture}
\end{center}
\vspace{-5mm}
\caption{An illustration of the factorisation $\al=\ob_P\be\ob_Q$ from Proposition \ref{prop:RPn_normalform}, where $\al\in\RP_{10}$.}
\label{fig:a=obo}
\end{figure}

\subsection{Semigroups and presentations}

We will be dealing extensively with both semigroup and monoid presentations, so we now take the time to fix our notation for these, as well as some general semigroup notions.  For further background on semigroups, the reader is referred to a monograph such as \cite{Hig} or \cite{Howie}.

An equivalence relation $\sim$ on a semigroup $S$ is a \emph{congruence} if $a\sim b$ and $c\sim d$ together imply $ac\sim bd$, for all $a,b,c,d\in S$.  If $\sim$ is a congruence on $S$, then the quotient $S/{\sim}$, which consists of all $\sim$-classes of $S$, is a semigroup under the natural induced operation.  The fundamental homomorphism theorem (for semigroups) states that if $\phi:S\to T$ is a semigroup homomorphism, then $S/\ker(\phi)\cong\im(\phi)$, where $\ker(\phi)$ is the congruence $\set{(a,b)\in S\times S}{a\phi=b\phi}$.

Let $X$ be an alphabet, and denote by $X^+$ (resp., $X^*$) the free semigroup (resp., free monoid) on~$X$.  If~$R\sub X^+\times X^+$ (resp., $R\sub X^*\times X^*$), we denote by $R^\sharp$ the congruence on $X^+$ (resp., $X^*$) generated by $R$.  We say a semigroup (resp., monoid) $S$ has \emph{semigroup} (resp., \emph{monoid}) \emph{presentation} $\pres{X\!}{\!R}$ if~${S\cong X^+/R^\sharp}$ (resp., $S\cong X^*/R^\sharp$) or, equivalently, if there is an epimorphism ${X^{+}}\to S$ (resp., $X^*\to S$) with kernel $R^\sharp$.  If $\phi$ is such an epimorphism, we say $S$ has \emph{presentation $\pres{X\!}{\!R}$ via $\phi$}.  A relation $(w_1,w_2)\in R$ will usually be displayed as an equation: $w_1=w_2$.  We will always be careful to specify whether a given presentation is a semigroup or monoid presentation.

We denote the \emph{empty word} (over any alphabet) by $1$ (so $X^*\sm X^+=\{1\}$ for any alphabet $X$).  If $w=x_1\cdots x_k$, where $x_1,\ldots,x_k\in X$, we write $\ell(w)=k$ for the \emph{length} of $w$.  
The word $x_i\cdots x_j$ is considered to be empty if either:
\bit
\item[(i)] $i>j$ and the subscripts are understood to be increasing; or
\item[(ii)] $i<j$ and the subscripts are understood to be decreasing.
\eit
When we are dealing with semigroup presentations, such a word will always be a subword of a larger (nonempty) word.

\subsection{Presentations for $\PnSn$ and $\P_n$}\label{sect:PnSn}

In Sections \ref{sect:RPnSn} and \ref{sect:RPn}, we will give presentations for $\RPnSn$ and $\RP_n$ (respectively).  To do this, it is crucial to know presentations for $\PnSn$ and $\P_n$, and we describe such presentations (from \cite{JEpnsn,JEgrpm}) in this section.
Consider alphabets $E=\{e_1,\ldots,e_n\}$ and $T=\set{t_{ij}}{\oijn}$, and define a (semigroup) homomorphism
\[
{\phi:(E\cup T)^+\to\PnSn}
\]
by
\[
e_i\phi =  \eb_i = \custpartn{1,3,4,5,7}{1,3,4,5,7}{\dotsups{1/3,5/7}\dotsdns{1/3,5/7}\stlines{1/1,3/3,5/5,7/7}\vertlabelshh{1/1,4/i,7/n}} \AND
t_{ij}\phi =  \tb_{ij} = \custpartn{1,3,4,5,7,8,9,11}{1,3,4,5,7,8,9,11}{\dotsups{1/3,5/7,9/11}\dotsdns{1/3,5/7,9/11}\stlines{1/1,3/3,4/4,5/5,7/7,8/8,9/9,11/11}\uarc48 \darc48\vertlabelshh{1/1,4/i,8/j,11/n}}.
\]
We will use symmetric notation when referring to the letters from $T$, so we write $t_{ij}=t_{ji}$ for all $\oijn$.
Consider the relations
\begin{align}
\tag{R1} e_i^2 &= e_i  &&\text{{for all $i$}}\\
\tag{R2} e_ie_j &= e_je_i   &&\text{{for distinct $i,j$}}\\
\tag{R3} t_{ij}^2 &= t_{ij}   &&\text{{for all $i,j$}}\\
\tag{R4} t_{ij}t_{kl} &= t_{kl}t_{ij}  &&\text{{for all $i,j,k,l$}}\\
\tag{R5} t_{ij}t_{jk} &= t_{jk}t_{ki}    &&\text{{for distinct $i,j,k$}}\\
\tag{R6} t_{ij}e_k &= e_kt_{ij}  &&\text{{if $k\not\in\{i,j\}$}}\\
\tag{R7} t_{ij}e_kt_{ij} &=  t_{ij} &&\text{{if $k\in\{i,j\}$}}\\
\tag{R8} e_kt_{ij}e_k &=  e_k &&\text{{if $k\in\{i,j\}$}}\\
\tag{R9} e_kt_{ki}e_it_{ij}e_jt_{jk}e_k &= e_kt_{kj}e_jt_{ji}e_it_{ik}e_k  &&\text{{for distinct $i,j,k$}}\\
\tag{R10} e_kt_{ki}e_it_{ij}e_jt_{jl}e_lt_{lk}e_k &= e_kt_{kl}e_lt_{li}e_it_{ij}e_jt_{jk}e_k &&\text{{for distinct $i,j,k,l$.}}
\end{align}
The following is \cite[Theorem 46]{JEpnsn}.

\ms
\begin{thm}\label{thm:PnSn}
The singular part of the partition monoid, $\PnSn$, has semigroup presentation
\[
\pres{E\cup T}{\text{\emph{(R1--R10)}}},
\]
via $\phi$. \epfres
\end{thm}

Now consider the alphabets $S=\{s_1,\ldots,s_{n-1}\}$, $E=\{e_1,\ldots,e_n\}$, $Q=\{q_1,\ldots,q_{n-1}\}$, and define a (monoid) homomorphism
\[
\psi:(S\cup E\cup Q)^*\to\P_n
\]
by
\[
s_i\psi=\sb_i = \custpartn{1,3,4,5,6,8}{1,3,4,5,6,8}{\dotsups{1/3,6/8}\dotsdns{1/3,6/8}\stlines{1/1,3/3,4/5,5/4,6/6,8/8}\vertlabelshh{1/1,4/i,8/n}}
\COMMA
e_i\psi=\eb_i = \custpartn{1,3,4,5,7}{1,3,4,5,7}{\dotsups{1/3,5/7}\dotsdns{1/3,5/7}\stlines{1/1,3/3,5/5,7/7}\vertlabelshh{1/1,4/i,7/n}}
\COMMA
q_i\psi=\qb_i = \tb_{i,i+1} = \custpartn{1,3,4,5,6,8}{1,3,4,5,6,8} {\dotsups{1/3,6/8}\dotsdns{1/3,6/8}\stlines{1/1,3/3,4/4,5/5,6/6,8/8}\uarc45\darc45\vertlabelshh{1/1,4/i,8/n}}
.
\]
Consider the relations
\begin{align}
\tag{R18} s_i^2 &= 1  &&\text{for all $i$}\\
\tag{R19} s_is_j &= s_js_i &&\text{if $|i-j|>1$}\\
\tag{R20} s_is_js_i &= s_js_is_j &&\text{if $|i-j|=1$}\\
\tag{R21} e_i^2 &= e_i  &&\text{for all $i$}   \\
\tag{R22} e_ie_j &= e_je_i   &&\text{for distinct $i,j$}    \\
\tag{R23} s_ie_j &= e_js_i &&\text{if $j\not\in\{i,i+1\}$}   \\
\tag{R24} s_ie_i &= e_{i+1}s_i  &&\text{for all $i$}   \\
\tag{R25} e_ie_{i+1}s_i &= e_ie_{i+1}  &&\text{for all $i$}   \\
\tag{R26}  q_i^2 &= q_i   &&\text{for all $i$}  \\
\tag{R27} q_iq_j &= q_jq_i   &&\text{for distinct $i,j$}   \\
\tag{R28} s_iq_j &= q_js_i  &&\text{if $|i-j|>1$}   \\
\tag{R29} s_is_jq_i &= q_js_is_j  &&\text{if $|i-j|=1$}   \\
\tag{R30} q_is_i = s_iq_i &= q_i  &&\text{for all $i$}   \\
\tag{R31} q_ie_j &= e_jq_i &&\text{if $j\not\in\{i,i+1\}$}   \\
\tag{R32} q_ie_jq_i &= q_i &&\text{if $j\in\{i,i+1\}$}  \\
\tag{R33} e_jq_ie_j &= e_j &&\text{if $j\in\{i,i+1\}$}.
\end{align}
(The jump in labels from relation (R10) to (R18) will become clear shortly.)
The next result was originally stated (in a different form) in \cite[Theorem 1.11]{HR2005}, and proved in \cite[Theorem 36]{JEgrpm}.

\ms
\begin{thm}\label{thm:Pn}
The partition monoid, $\P_n$, has monoid presentation 
\[
\pres{S\cup E\cup Q}{\text{\emph{(R18--R33)}}}
\]
via $\psi$. \epfres
\end{thm}

\section{The singular rook partition monoid $\RPnSn$}\label{sect:RPnSn}

In this section, we obtain a (semigroup) presentation for $\RPnSn$, the singular part of the rook partition monoid $\RP_n$; see Theorem \ref{thm:RPnSn}.  This presentation extends the presentation of $\PnSn$ stated in Theorem \ref{thm:PnSn}, and uses the minimum number of generators (Theorem \ref{thm:rankRPnSn}).  We note that the method used in this section could also be used to derive a presentation for $\RP_n$ itself (see Theorem \ref{thm:RPn}) from that of $\P_n$ (Theorem \ref{thm:Pn}).  
However, rather than duplicating the method for $\RP_n$, we instead use Theorem \ref{thm:RPnSn} as a stepping stone towards proving Theorem \ref{thm:RPn} in Section \ref{sect:RPn}, leading to a shorter proof.  It seems that this is the first time this method (of first finding a presentation for the singular subsemigroup) has been used, though we note some similarities to the author's recent work on (singular) symmetric inverse semigroups \cite{JEinsn2}.

Now define a new alphabet $O=\{o_1,\ldots,o_n\}$, and for each $i$, let
\[
\ob_i = \custpartn{1,3,5,7}{1,3,5,7}{\uvw4\lvw4\dotsups{1/3,5/7}\dotsdns{1/3,5/7}\stlines{1/1,3/3,5/5,7/7}\vertlabelshh{1/1,4/i,7/n}} \in\RP_n.
\]
Note that if $A=\{i_1,\ldots,i_k\}\sub\bn$, then $\ob_A$ (as defined in the previous section) may be factorised as $\ob_A=\ob_{i_1}\cdots\ob_{i_k}$; in particular, $\ob_i=\ob_{\{i\}}$ for any $i\in\bn$.  With this observation, the next result follows immediately from Proposition \ref{prop:RPn_normalform} and Theorem \ref{thm:PnSn}.

\ms
\begin{prop}\label{prop:RPnSn_gen}
The singular part of the rook partition monoid, $\RPnSn$, is generated (as a semigroup) by the set $\{\eb_1,\ldots,\eb_n\}\cup\set{\tb_{ij}}{\oijn}\cup\{\ob_1,\ldots,\ob_n\}$. \epfres
\end{prop}

\ms
\begin{rem}\label{rem:RPnSn_gen}
The previous result shows that $\RPnSn$ is generated by its idempotents; this property is shared by many \cite{BDP2002,JEpnsn,Maltcev2007,EastGray}, but not all \cite{DE2015,DEG2015,DE2016}, diagram monoids.  Under the above-mentioned embedding $\RP_n\to\P_{n+1}$, the generators $\ob_i\in\RP_n$ are mapped to $\tb_{i,n+1}\in\P_{n+1}$.  See also \cite{Howie1966}.
\end{rem}

By Proposition \ref{prop:RPnSn_gen}, we may define an epimorphism
\[
\Phi:(E\cup T\cup O)^+\to\RPnSn
\]
by $x\Phi=\xb$ for each $x\in E\cup T\cup O$.
Now consider the relations
\begin{align}
\tag{R11} o_i^2 &= o_i &&\text{for all $i$}   \\
\tag{R12} o_io_j &= o_jo_i  &&\text{for distinct $i,j$}    \\
\tag{R13} o_ie_j &= e_jo_i &&\text{for distinct $i,j$}   \\
\tag{R14} o_ie_io_i &= o_i &&\text{for all $i$}   \\
\tag{R15} e_io_ie_i &= e_i &&\text{for all $i$}   \\
\tag{R16} t_{ij}o_k &= o_kt_{ij} &&\text{for any $i,j,k$}   \\
\tag{R17} t_{ij}o_i=t_{ij}o_j &= o_io_j &&\text{for all $i,j$.}
\end{align}
Our aim in this section is to show that $\RPnSn$ has (semigroup) presentation $\pres{E\cup T\cup O}{\text{(R1--R17)}}$ via $\Phi$.

Since we already know $\Phi$ is surjective, it remains to show that $\ker\Phi$ is generated by the relations (R1--R17).  With this in mind, let $\sim$ be the congruence on $(E\cup T\cup O)^+$ generated by relations (R1--R17).
For $w\in(E\cup T\cup O)^+$, write $\wb=w\Phi\in\RPnSn$.  For convenience, we will also write $\overline{1}=1$ and $1\sim1$ (even though the empty word $1$ does not belong to $(E\cup T\cup O)^+$).

\ms
\begin{lemma}\label{lem:RPnSn_rels}
We have ${\sim}\sub\ker\Phi$.
\end{lemma}

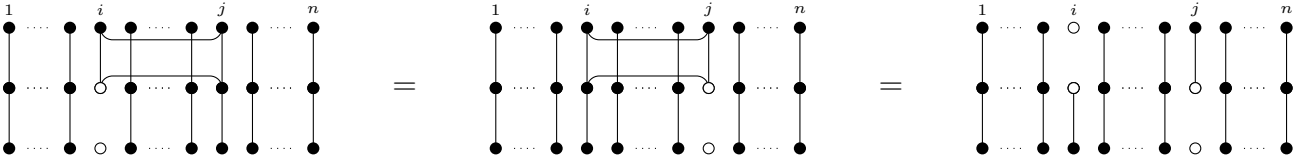
\begin{figure}[h]
\begin{center}
\begin{tikzpicture}[scale=.4]
\begin{scope}[shift={(0,2)}]	
\uverts{1,3,4,5,7,8,9,11}
\lverts{1,3,4,5,7,8,9,11}
\stlines{1/1,3/3,4/4,5/5,7/7,8/8,9/9,11/11}
\uarc48 
\darc48
\dotsups{1/3,5/7,9/11}
\dotsdns{1/3,5/7,9/11}
\vertlabelshh{1/1,4/i,8/j,11/n}
\end{scope}
\begin{scope}[shift={(0,0)}]	
\uverts{1,3,5,7,8,9,11}
\lverts{1,3,5,7,8,9,11}
\uvw4
\lvw4
\stlines{1/1,3/3,5/5,7/7,8/8,9/9,11/11}
\dotsups{1/3,5/7,9/11}
\dotsdns{1/3,5/7,9/11}
\draw(14,2)node{$=$};
\end{scope}
\begin{scope}[shift={(16,2)}]	
\uverts{1,3,4,5,7,8,9,11}
\lverts{1,3,4,5,7,8,9,11}
\stlines{1/1,3/3,4/4,5/5,7/7,8/8,9/9,11/11}
\uarc48 
\darc48
\dotsups{1/3,5/7,9/11}
\dotsdns{1/3,5/7,9/11}
\vertlabelshh{1/1,4/i,8/j,11/n}
\end{scope}
\begin{scope}[shift={(16,0)}]	
\uverts{1,3,4,5,7,9,11}
\lverts{1,3,4,5,7,9,11}
\uvw8
\lvw8
\stlines{1/1,3/3,4/4,5/5,7/7,9/9,11/11}
\dotsups{1/3,5/7,9/11}
\dotsdns{1/3,5/7,9/11}
\draw(14,2)node{$=$};
\end{scope}
\begin{scope}[shift={(32,2)}]	
\vertlabelshh{1/1,4/i,8/j,11/n}
\uverts{1,3,5,7,8,9,11}
\lverts{1,3,5,7,8,9,11}
\uvw4
\lvw4
\stlines{1/1,3/3,5/5,7/7,8/8,9/9,11/11}
\dotsups{1/3,5/7,9/11}
\dotsdns{1/3,5/7,9/11}
\end{scope}
\begin{scope}[shift={(32,0)}]	
\uverts{1,3,4,5,7,9,11}
\lverts{1,3,4,5,7,9,11}
\uvw8
\lvw8
\stlines{1/1,3/3,4/4,5/5,7/7,9/9,11/11}
\uvw4
\dotsups{1/3,5/7,9/11}
\dotsdns{1/3,5/7,9/11}
\end{scope}
\end{tikzpicture}
\end{center}
\vspace{-5mm}
\caption{Diagrammatic proof of relation (R17): $\tb_{ij}\ob_i=\tb_{ij}\ob_j=\ob_i\ob_j$.}
\label{fig:R17}
\end{figure}

\pf We need to show that each of relations (R1--R17) hold as equations in $\RPnSn$ when the words are replaced by their images under $\Phi$.  We already know from Theorem \ref{thm:PnSn} that this is the case for (R1--R10).  The remaining relations may easily be checked diagramatically; we do this for (R17) in Figure \ref{fig:R17}, and leave the rest for the reader. \epf

Establishing the reverse inclusion, $\ker\Phi\sub{\sim}$, forms the bulk of this section.  The main step in doing this is to obtain a ``word version'' of Proposition \ref{prop:RPn_normalform}; see Proposition \ref{prop:AuB}.  Our first aim is to show that any word over $(E\cup T\cup O)^+$ is $\sim$-equivalent to an element of $O^*(E\cup T)^*O^*$; see Lemma~\ref{lem:w123}, the proof of which requires the next technical result.

\ms
\begin{lemma}\label{lem:woe}
Let $w\in (E\cup T)^*$ and $i\in\bn$.  Then $wo_ie_i\sim w_1w_2w_3$ for some $w_1,w_3\in O^*$ and $w_2\in (E\cup T)^*$.
\end{lemma}

\pf We prove the result by induction on $\ell(w)$, the length of $w$.  If $\ell(w)=0$, then we are done (with $w_1=o_i$, $w_2=e_i$ and $w_3=1$), so suppose $\ell(w)\geq1$, and write $w=ux$, where $x\in E\cup T$ (so $u\in(E\cup T)^*$ and $\ell(u)=\ell(w)-1$).  We now consider separate cases, according to whether $x$ belongs to $E$ or $T$.

{\bf Case 1.}  First suppose $x=e_j\in E$.  If $j=i$, then $wo_ie_i=ue_io_ie_i\sim ue_i=w$, by (R15), and we are done (with $w_1=w_3=1$ and $w_2=w$).  If $j\not=i$, then $wo_ie_i=ue_jo_ie_i \sim uo_ie_ie_j$, by (R2) and (R13), and we are done after applying an induction hypothesis to $uo_ie_i$.

{\bf Case 2.}  Next suppose $x=t_{jk}\in T$.  If $i\not\in\{j,k\}$, then $wo_ie_i=ut_{jk}o_ie_i\sim uo_ie_it_{jk}$, by (R6) and (R16), and again we are done after applying an induction hypothesis.  If $i\in\{j,k\}$, then, writing $\{l\}=\{j,k\}\sm\{i\}$, we have $wo_ie_i=ut_{jk}o_ie_i\sim ut_{jk}o_le_i\sim ut_{jk}e_io_l=we_io_l$, by (R17) and (R13), and we are done (with $w_1=1$, $w_2=we_i$ and $w_3=o_l$).  \epf

\ms
\begin{lemma}\label{lem:w123}
Let $w\in(E\cup T\cup O)^+$.  Then $w\sim w_1w_2w_3$ for some $w_1,w_3\in O^*$ and $w_2\in (E\cup T)^*$.
\end{lemma}

\pf This is clearly true if $\ell(w)=1$, so suppose $\ell(w)\geq2$ and write $w=ux$, where $x\in E\cup T\cup O$.  By an induction hypothesis, $u\sim u_1u_2u_3$ for some $u_1,u_3\in O^*$ and $u_2\in (E\cup T)^*$.  

{\bf Case 1.}  If $x\in O$, then we are done, with $w_1=u_1$, $w_2=u_2$ and $w_3=u_3x$.

{\bf Case 2.}  If $x\in T$, then relation (R16) gives $u_3x\sim xu_3$, and we are done, with $w_1=u_1$, $w_2=u_2x$ and $w_3=u_3$.

{\bf Case 3.}  Finally, suppose $x=e_i\in E$.  If $\ell(u_3)=0$, then we are done, so suppose $\ell(u_3)\geq1$, and write $u_3=o_{j_1}\cdots o_{j_k}$; by (R11) and (R12), we may assume that $j_1<\cdots<j_k$.  If $i\not\in\{j_1,\ldots,j_k\}$, then $u_3e_i\sim e_iu_3$, by (R13), and we are done.
So suppose $i\in\{j_1,\ldots,j_k\}$; say, $i=j_l$.  Then $u_3e_i\sim o_ie_iv$, by (R12) and (R13), where $v=o_{j_1}\cdots o_{j_{l-1}}o_{j_{l+1}}\cdots o_{j_k}$.  By Lemma \ref{lem:woe}, $u_2o_ie_i\sim v_1v_2v_3$ for some $v_1,v_3\in O^*$ and $v_2\in (E\cup T)^*$.  Putting all this together, we have $w\sim u_1v_1v_2v_3v$, and we are done, with $w_1=u_1v_1$, $w_2=v_2$ and $w_3=v_3v$. \epf

Next we wish to show (in Proposition \ref{prop:AuB}) that the words $w_1,w_2,w_3$ from Lemma \ref{lem:w123} may be chosen so that $\wb_1,\wb_2,\wb_3$ correspond to the factorisation of $\wb\in\RPnSn$ as in Proposition \ref{prop:RPn_normalform}.  To achieve this goal, we need several intermediate results.

For $A=\{i_1,\ldots,i_k\}\sub\bn$ with $i_1<\cdots<i_k$, we define the words
\[
o_A=o_{i_1}\cdots o_{i_k}\in O^* \AND t_A=t_{i_1i_2}t_{i_2i_3}\cdots t_{i_{k-1}i_k}\in T^*.
\]
Note that $o_A=1$ if $A=\emptyset$, while $t_A=1$ if $|A|\leq1$.
Note also that $o_A\Phi=\ob_A$, agreeing with our earlier use of this notation.

\ms
\begin{lemma}\label{lem:ot}
Let $A\sub\bn$ and let $i\in A$.  Then $o_A\sim o_it_A\sim t_Ao_i$.
\end{lemma}

\pf By (R16), it suffices to show that $o_A\sim o_it_A$, and we do this by induction on $|A|$.  If $|A|=1$, then there is nothing to show, since then $o_A=o_i$ and $t_A=1$.  So suppose $|A|\geq2$, and write $A=\{j_1,\ldots,j_k\}$ with $j_1<\cdots<j_k$.  Then $i=j_l$ for some $1\leq l\leq k$, and we write $B=A\sm\{i\}$.  Since $k=|A|\geq2$, we must have $l>1$ or $l<k$ (or both).  We will assume that $l<k$ (the other case is similar).  Then
\begin{align*}
o_A &\sim o_io_B &&\text{by (R12)}\\
&\sim o_io_{j_{l+1}}t_B &&\text{by an induction hypothesis}\\
&\sim o_it_{ij_{l+1}}(t_{j_1j_2}\cdots t_{j_{l-2}j_{l-1}} t_{j_{l-1}j_{l+1}} t_{j_{l+1}j_{l+2}}\cdots t_{j_{k-1}j_k}) &&\text{by (R16) and (R17)}\\
&\sim o_i(t_{j_1j_2}\cdots t_{j_{l-2}j_{l-1}}t_{ij_{l+1}}t_{j_{l-1}j_{l+1}} t_{j_{l+1}j_{l+2}}\cdots t_{j_{k-1}j_k}) &&\text{by (R4)}\\
&= o_i(t_{j_1j_2}\cdots t_{j_{l-2}j_{l-1}}t_{j_lj_{l+1}}t_{j_{l+1}j_{l-1}} t_{j_{l+1}j_{l+2}}\cdots t_{j_{k-1}j_k}) \\ 
&\sim o_i(t_{j_1j_2}\cdots t_{j_{l-2}j_{l-1}}t_{j_{l-1}j_l}t_{j_lj_{l+1}} t_{j_{l+1}j_{l+2}}\cdots t_{j_{k-1}j_k}) &&\text{by (R5)}\\
\epfreseq &= o_it_A.
\end{align*}

Using the previous result, we may now strengthen Lemma \ref{lem:w123}.

\ms
\begin{cor}\label{cor:w123}
Let $w\in(E\cup T\cup O)^+$.  Then $w\sim w_1w_2w_3$ for some $w_1,w_3\in O^*$ and $w_2\in (E\cup T)^*$, with $\ell(w_1),\ell(w_3)\leq1$.
\end{cor}

\pf By Lemma \ref{lem:w123}, $w\sim u_1u_2u_3$ for some $u_1,u_3\in O^*$ and $u_2\in (E\cup T)^*$.  Then $u_1\sim o_A$ and $u_3\sim o_B$ for some (possibly empty) subsets $A,B\sub\bn$, by (R11) and (R12).  If $A\not=\emptyset$, then choose some $i\in A$; if $B\not=\emptyset$, then choose some $j\in B$.  Then 
\[
w\sim u_1u_2u_3\sim o_Au_2o_B \sim 
\begin{cases}
u_2 &\text{if $A=\emptyset=B$}\\
o_i(t_Au_2) &\text{if $A\not=\emptyset=B$}\\
(u_2t_B)o_j &\text{if $A=\emptyset\not=B$}\\
o_i(t_Au_2t_B)o_j &\text{if $A\not=\emptyset\not=B$,}
\end{cases}
\]
by Lemma \ref{lem:ot}, and the proof is complete. \epf

Next we wish to show how the relations may be used to move a generator from $O$ through a word from $(E\cup T)^*$ in certain circumstances.  To do this, we require a result from \cite{JEinsn2} on (singular) symmetric inverse semigroups.

For distinct $i,j\in\bn$, we define the word $e_{ij}=e_it_{ij}e_j\in(E\cup T)^+$, recalling that we are using symmetric notation for $t_{ij}=t_{ji}$.  One easily checks that
\[
\eb_{ij} = \begin{cases}
\custpartn{1,3,4,5,7,8,9,11}{1,3,4,5,7,8,9,11}{\dotsups{1/3,5/7,9/11}\dotsdns{1/3,5/7,9/11}\stlines{1/1,3/3,5/5,7/7,9/9,11/11,8/4}\vertlabelshh{1/1,4/i,8/j,11/n}} &\text{if $i<j$}\\~\\
\custpartn{1,3,4,5,7,8,9,11}{1,3,4,5,7,8,9,11}{\dotsups{1/3,5/7,9/11}\dotsdns{1/3,5/7,9/11}\stlines{1/1,3/3,5/5,7/7,9/9,11/11,4/8}\vertlabelshh{1/1,4/j,8/i,11/n}} &\text{if $j<i$.}
\end{cases}
\]
In particular, $\eb_{ij}\not=\eb_{ji}$.  The next result is \cite[Proposition 2.2]{JEinsn2}.

\ms
\begin{prop}\label{prop:InSn}
The singular part of the symmetric inverse monoid, $\InSn$, is generated (as a semigroup) by the set $\set{\eb_{ij}}{i,j\in\bn,\ i\not=j}$. \epfres
\end{prop}

\ms
\begin{rem}
Defining relations were also given in \cite[Theorem 2.1]{JEinsn2} but we do not need those here.
\end{rem}

\ms
\begin{lemma}\label{lem:okeij}
Let $i,j\in\bn$ with $i\not=j$, and let $k\in\{i\}^c$.  Then 
\[
o_ke_{ij} \sim \begin{cases}
e_{ij}o_k &\text{if $k\not=j$}\\
e_{ij}o_i &\text{if $k=j$.}
\end{cases}
\]
\end{lemma}

\pf The case in which $k\not=j$ follows immediately from (R13) and (R16).  We also have 
\[
o_je_{ij} = o_j e_i t_{ij} e_j \sim e_i o_j t_{ij} e_j \sim e_i t_{ij} o_j e_j \sim e_i t_{ij} o_i e_j \sim e_i t_{ij} e_j o_i = e_{ij}o_i,
\]
by (R13), (R16), (R17), (R13), respectively. \epf

As usual, for $\al\in\I_n$ and $i\in\dom(\al)$, we write $i\al$ for the unique element of $\codom(\al)$ such that $\{i,(i\al)'\}$ is a block of $\al$.
Note that Lemma \ref{lem:okeij} says that $o_ke_{ij}\sim e_{ij}o_{k\eb_{ij}}$ for all $k\in\{i\}^c=\dom(\eb_{ij})$; compare with the above illustration(s) of $\eb_{ij}$.

\ms
\begin{cor}\label{cor:oiw}
Let $w\in(E\cup T)^*$ be such that $\wb\in\I_n$.  Then $o_iw\sim wo_{i\wb}$ for any $i\in\dom(\wb)$.
\end{cor}

\pf The result is trivial if $\ell(w)=0$, so suppose $\ell(w)\geq1$.  In particular, $\wb\in\InSn$, so $\wb=\eb_{i_1j_1}\cdots\eb_{i_kj_k}$ for some $i_t,j_t$, by Proposition \ref{prop:InSn}.  It follows from Theorem \ref{thm:PnSn} that $w\sim e_{i_1j_1}\cdots e_{i_kj_k}$, and the result now follows from Lemma \ref{lem:okeij} and a simple induction on $k$. \epf

For the proof of the next result, note that if
\[
\al=\partnlistnodefshort{A_1}{A_r}{A_1}{A_r}\in\J_n,
\]
then $\al=\overline{t_{A_1}\cdots t_{A_r}}$, where the words $t_A$ were defined before Lemma \ref{lem:ot}.

\ms
\begin{lemma}\label{lem:owwo}
Suppose $i,j\in\bn$ and $w\in (E\cup T)^*$ are such that $i$ and $j'$ belong to the same block $A\cup B'$ of $\wb$.  Then $o_iw\sim wo_j\sim o_iuo_j$ for some $u\in(E\cup T)^+$  such that $\ub = (\wb\sm\{A\cup B'\})\cup\{A,B'\}$.
\end{lemma}

\pf Suppose first that $\ell(w)=0$.  It follows that $w=1$ and $i=j$, in which case $o_iw= wo_j=o_i\sim o_ie_io_i$, by~(R14), and we are done (with $u=e_i$).  For the remainder of the proof, we assume that $\ell(w)\geq1$.  In particular, $\wb\in\PnSn$.  Write
\[
\wb = \partnlist{A_1}{A_r}{C_1}{C_p}{B_1}{B_r}{D_1}{D_q},
\]
and let $\be,\ga,\de$ be such that $\wb=\be\ga\de$, as in Proposition \ref{prop:Pn_normalform}.  Without loss of generality, we may suppose that $A=A_r$, $B=B_r$, and that $a_r=i$ and $b_r=j$.  
Let $w_1,w_3\in T^*$ and $w_2\in (E\cup T)^+$ be such that $\wb_1=\be$, $\wb_2=\ga$ and $\wb_3=\de$.  Note that $w_2\not=1$ since $r=\rank(\wb)<n$, and that $i\wb_2=i\ga=j$.  Then
\begin{align*}
o_iw &\sim o_i w_1w_2w_3 &&\text{by Theorem \ref{thm:PnSn}}\\
&\sim w_1o_i w_2w_3 &&\text{by (R16)}\\
&\sim w_1 w_2o_jw_3 &&\text{by Corollary \ref{cor:oiw}}\\
&\sim w_1 w_2w_3o_j &&\text{by (R16)}\\
&\sim wo_j &&\text{by Theorem \ref{thm:PnSn}.}
\end{align*}
\newpage
We also have
\begin{align*}
o_iw &\sim w_1o_i w_2w_3 &&\text{by Theorem \ref{thm:PnSn} and (R16), as above}\\
&\sim w_1o_ie_io_i w_2w_3 &&\text{by (R14)}\\
&\sim o_iw_1e_i w_2w_3o_j &&\text{by Corollary \ref{cor:oiw} and (R16), as above}\\
&=o_iuo_j,
\end{align*}
where $u=w_1e_i w_2w_3$.  Note that
\[
\overline{e_iw_2} = \partnlistnodefshort{a_1}{a_{r-1}}{b_1}{b_{r-1}} \qquad\text{so that}\qquad \ub=\overline{w_1e_i w_2w_3}= \left( \begin{array}{c|c|c|c|c|c|c} \!\! A_1 & \cdots & A_{r-1} & A & C_1 & \cdots & C_p\ \ \\ \cline{4-7}
\!\! B_1 & \cdots & B_{r-1} & B & D_1 & \cdots & D_q \ \
\end{array} \!\!\! \right),
\]
completing the proof. \epf

We are now ready to describe the promised normal forms for words over $E\cup T\cup O$.
As in Proposition~\ref{prop:RPn_normalform}, for the statement of the next result, we will allow ourselves to refer to the empty set $\emptyset$ as a block of any element of $\P_n$.

\ms
\begin{prop}\label{prop:AuB}
Let $w\in(E\cup T\cup O)^+$.  Then $w\sim o_Auo_B$ for some (possibly empty) $A,B\sub\bn$ and $u\in(E\cup T)^+$ with $A$ and $B'$ (possibly empty) blocks of $\ub$.
\end{prop}

\pf By Corollary \ref{cor:w123}, $w\sim w_1vw_2$ for some $w_1,w_2\in O^*$ and $v\in (E\cup T)^*$, with $\ell(w_1),\ell(w_2)\leq1$.  If $\ell(w_1)=\ell(w_2)=0$, we are already done (with $A=B=\emptyset$ and $u=w$), so suppose this is not the case.  

{\bf Case 1.}  Suppose first that $\ell(w_1)=1$ and $\ell(w_2)=0$, so $w\sim o_iv$ for some $i\in\bn$.  Write $[i]_{\vb}=A\cup B'$, where $A,B\sub\bn$, noting that $A$ is non-empty (since $i\in A$), but $B$ may be empty.   Since $A\cup B'$ is a block of $\vb$, it follows that $\tb_A\vb\tb_B=\vb$, so Theorem \ref{thm:PnSn} gives $v\sim t_Avt_B$ (but note that this just says~$1\sim1$ if $\ell(v)=0$).

{\bf Subcase 1.1.}  Suppose first that $B=\emptyset$.  Then by Lemma \ref{lem:ot} and the above calculations, we have
$
w\sim o_iv \sim o_it_Av \sim o_Av,
$
so the proof is complete in this case (with $u=v$).

{\bf Subcase 1.2.}  Now suppose $B\not=\emptyset$, and choose some $j\in B$.  By Lemma \ref{lem:owwo}, $w\sim o_iv\sim o_iuo_j$ for some $u\in(E\cup T)^+$ such that $\ub = (\vb\sm\{A\cup B'\})\cup\{A,B'\}$.  Since $A$ and $B'$ are blocks of $\ub$, it again follows from Theorem \ref{thm:PnSn} that $u\sim t_Aut_B$.  Together with Lemma \ref{lem:ot}, we then obtain
$
w\sim o_iuo_j\sim o_it_Aut_Bo_j\sim o_Auo_B,
$
completing the proof in this case.

{\bf Case 2.}  The case in which $\ell(w_1)=0$ and $\ell(w_2)=1$ is similar to the previous case.

{\bf Case 3.}  Finally, suppose $\ell(w_1)=\ell(w_2)=1$, so $w\sim o_ivo_j$ for some $i,j\in\bn$.  Write $[i]_{\vb}=A\cup C'$ and $[j']_{\vb}=D\cup B'$, noting that $A$ and $B$ are non-empty, but that $C$ and/or $D$ may be empty.  (Note that it is possible that $[i]_{\vb}=[j']_{\vb}$, in which case $A=D$ and $B=C$, but this is included in the case that $C$ and $D$ are non-empty.)

{\bf Subcase 3.1.}  Suppose first that $D\not=\emptyset$, and let $k\in D$.  Then $vo_j\sim o_kv$, by Lemma \ref{lem:owwo}, so
\[
w\sim o_ivo_j\sim o_io_kv \sim \begin{cases}
o_iv &\text{by (R11), if $k=i$}\\
o_it_{ik}v &\text{by (R17) and (R16), if $k\not=i$.}
\end{cases}
\]
In either case, this reduces to Case 1.

{\bf Subcase 3.2.}  In similar fashion to the previous subcase, if $C\not=\emptyset$, then we may reduce to Case 2.

{\bf Subcase 3.3.}  Finally, suppose $C=D=\emptyset$.  Then
$w\sim o_ivo_j \sim o_it_Avt_Bo_j \sim o_Avo_B$,
by Theorem \ref{thm:PnSn} and Lemma \ref{lem:ot}, and we are done, with $u=v$. \epf

We may now prove the main result of this section.

\ms
\begin{thm}\label{thm:RPnSn}
The singular part of the rook partition monoid, $\RPnSn$, has semigroup presentation
\[
\pres{E\cup T\cup O}{\text{\emph{(R1--R17)}}}
\]
via $\Phi$.
\end{thm}

\pf It remains to show that $\ker\Phi\sub{\sim}$, so suppose $w_1,w_2\in(E\cup T\cup O)^+$ are such that $\wb_1=\wb_2$.  By Proposition \ref{prop:AuB}, $w_1\sim o_Au_1o_B$ and $w_2\sim o_Cu_2o_D$, for some (possibly empty) $A,B,C,D\sub\bn$ and $u_1,u_2\in(E\cup T)^+$ such that $A$ and $B'$ are (possibly empty) blocks of $\ub_1$, and $C$ and $D'$ are (possibly empty) blocks of $\ub_2$.  But then $\ob_A\ub_1\ob_B = \wb_1 = \wb_2 = \ob_C\ub_2\ob_D$.  Proposition \ref{prop:RPn_normalform} then gives $A=C$, $B=D$ and $\ub_1=\ub_2$.  Theorem \ref{thm:PnSn} gives $u_1\sim u_2$.  Putting this all together, we have $w_1 \sim o_Au_1o_B \sim o_Au_2o_B = o_Cu_2o_D \sim w_2$, completing the proof. \epf

The next result shows that the generating set used in Theorem \ref{thm:RPnSn} has the smallest possible size among all generating sets for $\RPnSn$.  Recall that the \emph{rank} a semigroup $S$, denoted $\rank(S)$, is the smallest size of a generating set for $S$.  If $S$ is generated by its idempotents, then the \emph{idempotent rank} of $S$, denoted $\idrank(S)$, is the smallest size of a generating set consisting entirely of idempotents.  For the proof of the next result, for $\oijn$, let $\ve_{ij}$ be the equivalence relation on $\bn$ with~$\{i,j\}$ as its only non-trivial equivalence class.

\ms
\begin{thm}\label{thm:rankRPnSn}
We have $\ds{\rank(\RPnSn)=\idrank(\RPnSn)=\frac{n^2+3n}2}$. \epfres
\end{thm}

\pf We know that $(E\cup T\cup O)\Phi$ is an idempotent generating set for $\RPnSn$.  Since $|E\cup T\cup O|=n+\binom{n}{2}+n=\frac{n^2+3n}2$, it follows that $\idrank(\RPnSn)\leq\frac{n^2+3n}2$.  Since also $\rank(\RPnSn)\leq\idrank(\RPnSn)$, it remains to show that any generating set for $\RPnSn$ has size at least $\frac{n^2+3n}2$.  With this in mind, suppose $\RPnSn=\la\Si\ra$.  We claim that:
\bit
\item[(i)] for all $i\in\bn$, $\Si$ contains an element with domain $\{i\}^c$ and the block $\{i\}$;
\item[(ii)] for all $i\in\bn$, $\Si$ contains an element with domain $\{i\}^c$ and the rook dot $i$; and
\item[(iii)] for all $\oijn$, $\Si$ contains an element with domain $\bn$ and kernel $\ve_{ij}$.
\eit
In fact, the proofs of these are all very similar, so we just prove (i).  With this in mind, let ${i\in\bn}$, and consider an expression $\eb_i=\al_1\cdots\al_k$, where $\al_1,\ldots,\al_k\in\Si$.  Then
$
n-1=\rank(\eb_i)=\rank(\al_1\cdots\al_k)\leq\rank(\al_1)\leq n-1,
$
so that $\rank(\al_1)=n-1$.  We also have
$
\{i\}^c=\dom(\eb_i)=\dom(\al_1\cdots\al_k)\sub\dom(\al_1),
$
so that $\dom(\al_1)$ is either $\{i\}^c$ or $\bn$.  If the latter was the case, then (together with $\rank(\al_1)=n-1$) this would imply that $\ker(\al_1)=\ve_{uv}$ for some $1\leq u<v\leq n$, and this would give
$
\ve_{uv}=\ker(\al_1)\sub\ker(\al_1\cdots\al_k)=\ker(\eb_i)=\De,
$
a contradiction.
So, in fact, $\dom(\al_1)=\{i\}^c$.  It follows that either $\{i\}$ is a block of $\al_1$ or else $i$ is a rook dot of $\al_1$.  If the latter was the case, then $i$ would also be a rook dot of $\al_1\cdots\al_k=\eb_i$, a contradiction.  So, in fact, $\{i\}$ is a block of $\al_1$.  This completes the proof of (i).  Finally, we note that $\Si$ has (at least) $n$ elements of type (i), $n$ of type (ii), and $\binom{n}{2}$ of type (iii), so that $|\Si|\geq n+n+\binom{n}{2}=\frac{n^2+3n}2$.  
\epf

\ms
\begin{rem}
Theorem \ref{thm:rankRPnSn} may also be proved using the general result \cite[Theorem 5.10]{EastGray}, utilising \emph{Green's relations} and the \emph{regular $*$-semigroup} structure of $\RP_n$, which we have not investigated here.  It would be interesting to study $\RP_n$ with the machinery developed in \cite{EastGray}, as it would allow one to study \emph{all} the ideals of $\RP_n$ (with $\RPnSn$ being the largest proper ideal).  In particular, the following questions seem worthy of study:
\bit
\item[(i)] What are the ranks of the other proper ideals of $\RP_n$?
\item[(ii)] Are the other proper ideals of $\RP_n$ generated by their idempotents?  If so, what are their idempotent ranks?  Are the ranks and idempotent ranks equal?
\item[(iii)] How many minimal-size (idempotent) generating sets are there for $\RPnSn$?
\eit
These questions are all answered for the proper ideals of partition monoid $\P_n$ in \cite{EastGray}, as well as for the \emph{Brauer monoid} and \emph{Jones monoid}.  See also \cite{Howie1978}.
\end{rem}

\section{The rook partition monoid $\RP_n$}\label{sect:RPn}

The goal of this section is to obtain a (monoid) presentation for the rook partition monoid $\RP_n$; in fact, we obtain several such presentations (see Theorems \ref{thm:RPn}, \ref{thm:RPn2} and \ref{thm:RPn3}).  
Our approach makes crucial use of Theorem \ref{thm:RPnSn}.  Specifically, we first postulate a presentation for $\RP_n$ (built up from the presentation of $\P_n$ in Theorem \ref{thm:Pn}), and then show that the presentation for $\RPnSn$ from Theorem~\ref{thm:RPnSn} can be \emph{embedded} (in a certain sense) in the stated presentation for $\RP_n$ (see Lemmas~\ref{lem:reduction1} and \ref{lem:reduction2}).

As with Proposition \ref{prop:RPnSn_gen}, the next result follows from Proposition \ref{prop:RPn_normalform} and Theorem~\ref{thm:Pn}.

\ms
\begin{prop}\label{prop:RPn_gen}
The rook partition monoid, $\RP_n$, is generated (as a monoid) by the set
\[\epfreseq
\{\sb_1,\ldots,\sb_{n-1}\}\cup \{\eb_1,\ldots,\eb_n\}\cup\{\qb_1,\ldots,\qb_{n-1}\}\cup\{\ob_1,\ldots,\ob_n\}.
\]
\end{prop}

Consider the alphabets
\[
S=\{s_1,\ldots,s_{n-1}\} \COMMA E=\{e_1,\ldots,e_n\} \COMMA Q=\{q_1,\ldots,q_{n-1}\}\COMMA O=\{o_1,\ldots,o_n\} ,
\]
as defined earlier.  By Proposition \ref{prop:RPn_gen}, we may define a (monoid) epimorphism
\[
\Psi:(S\cup E\cup Q\cup O)^*\to\RP_n
\]
by $x\Psi=\xb$ for each $x\in S\cup E\cup Q\cup O$.
Consider also the relations
\begin{align*}
\tag{R34} o_i^2 &= o_i  &&\text{for all $i$}   \\
\tag{R35} o_io_j &= o_jo_i   &&\text{for distinct $i,j$}    \\
\tag{R36} s_io_j &= o_js_i &&\text{if $j\not\in\{i,i+1\}$}   \\
\tag{R37} s_io_i &= o_{i+1}s_i  &&\text{for all $i$}   \\
\tag{R38} o_io_{i+1}s_i &= o_io_{i+1}  &&\text{for all $i$}   \\
\tag{R39} o_ie_j &= e_jo_i &&\text{for distinct $i,j$}   \\
\tag{R40} o_ie_io_i &= o_i  &&\text{for all $i$}   \\
\tag{R41} e_io_ie_i &= e_i  &&\text{for all $i$}   \\
\tag{R42} q_io_j &= o_jq_i &&\text{for all $i,j$}   \\
\tag{R43} q_io_i = q_io_{i+1} &= o_io_{i+1} &&\text{for all $i$.}
\end{align*}
Our goal in this section is to show that $\RP_n$ has (monoid) presentation $\pres{S\cup E\cup Q\cup O}{\text{(R18--R43)}}$ via $\Psi$.

Since we already know $\Psi$ is surjective, it remains to show that $\ker\Psi$ is generated by the relations (R18--R43).  With this in mind, let $\approx$ denote the congruence on $(S\cup E\cup Q\cup O)^*$ generated by relations (R18--R43).  Without causing confusion, we will write $\wb=w\Psi$ for all $w\in(S\cup E\cup Q\cup O)^*$.
As with Lemma \ref{lem:RPnSn_rels}, the next result may easily be proved diagrammatically.

\ms
\begin{lemma}\label{lem:RPn_rels}
We have ${\approx}\sub\ker\Psi$. \epfres
\end{lemma}

To establish the reverse containment, $\ker\Psi\sub{\approx}$, we need to show that $\ub=\vb\implies u\approx v$ for all $u,v\in(S\cup E\cup Q\cup O)^*$.  The next result shows that this is the case for words $u,v$ in certain subsets of $(S\cup E\cup Q\cup O)^*$.

\ms
\begin{prop}\label{prop:InPn}
If $u,v\in(S\cup E\cup Q)^*$ or $u,v\in(S\cup O)^*$, then $\ub=\vb\implies u\approx v$.
\end{prop}

\pf This is clearly the case for $u,v\in(S\cup E\cup Q)^*$, since (R18--R33) are defining relations for~$\P_n=\la (S\cup E\cup Q)\Psi\ra$ (Theorem \ref{thm:Pn}).  Similarly, (R18--R20, R34--R38) constitute defining relations for $\R_n=\la (S\cup O)\Psi\ra$ (see \cite[Theorem 4.8]{Gilbert2006} and \cite[Proposition 3.22]{JEboppp}). \epf

\begin{figure}[h]
\begin{center}
\begin{tikzpicture}[scale=.4]
\begin{scope}[shift={(0,2)}]	
\uverts{1,3,4,5,7,8,9,11}
\lverts{1,3,4,5,6,8,9,11}
\stlines{1/1,3/3,4/4,5/6,7/8,8/5,9/9,11/11}
\dotsups{1/3,5/7,9/11}
\dotsdns{1/3,6/8,9/11}
\draw(-.2,1)node[left]{$\overline{\si}_{ij}^{-1}$};
\draw[|-|] (-.2,2)--(-.2,0);
\vertlabelshh{1/1,4/i,8/j,11/n}
\end{scope}
\begin{scope}[shift={(0,0)}]	
\uverts{1,3,4,5,6,8,9,11}
\lverts{1,3,4,5,6,8,9,11}
\uarc45
\darc45
\stlines{1/1,3/3,4/4,5/5,6/6,8/8,9/9,11/11}
\dotsups{1/3,6/8,9/11}
\dotsdns{1/3,6/8,9/11}
\draw(-.2,1)node[left]{$\qb_i$};
\draw (-.2,2)--(-.2,0);
\draw(14,1)node{$=$};
\end{scope}
\begin{scope}[shift={(0,-2)}]	
\lverts{1,3,4,5,7,8,9,11}
\uverts{1,3,4,5,6,8,9,11}
\stlines{1/1,3/3,4/4,6/5,8/7,5/8,9/9,11/11}
\dotsdns{1/3,5/7,9/11}
\dotsups{1/3,6/8,9/11}
\draw(-.2,1)node[left]{$\overline{\si}_{ij}$};
\draw[|-|] (-.2,2)--(-.2,0);
\end{scope}
\begin{scope}[shift={(16,0)}]	
\uverts{1,3,4,5,7,8,9,11}
\lverts{1,3,4,5,7,8,9,11}
\stlines{1/1,3/3,4/4,5/5,7/7,8/8,9/9,11/11}
\uarc48 
\darc48
\dotsups{1/3,5/7,9/11}
\dotsdns{1/3,5/7,9/11}
\draw(12.2,1)node[right]{$\tb_{ij}$};
\draw[|-|] (12.2,2)--(12.2,0);
\vertlabelshh{1/1,4/i,8/j,11/n}
\end{scope}
\end{tikzpicture}
\end{center}
\vspace{-5mm}
\caption{Diagrammatic proof that $\overline{\tau}_{ij}=\tb_{ij}$.}
\label{fig:tauij}
\end{figure}
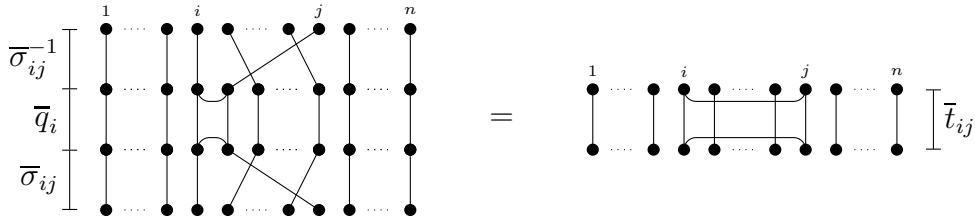

For a word $w=s_{i_1}\cdots s_{i_k}\in S^*$, we write $w^{-1}=s_{i_k}\cdots s_{i_1}$.  Note that $ww^{-1}\approx w^{-1}w\approx 1$ for all $w\in S^*$, by (R18).
For $\oijn$, define words $\si_{ij}=s_{i+1}\cdots s_{j-1}$ (which is empty if $j=i+1$) and $\tau_{ij}=\tau_{ji}=\si_{ij}^{-1}q_i\si_{ij}$.  Note that $\overline{\tau}_{ij}=\tb_{ij}$, as shown in Figure \ref{fig:tauij}.  Note also that $\tau_{i,i+1}=q_i$ for all~$i$.

As stated earlier, our strategy involves (somehow) linking the presentations $\pres{E\cup T\cup O}{\text{(R1--R17)}}$ and $\pres{S\cup E\cup Q\cup O}{\text{(R18--R43)}}$.  In order to make this link explicit, we define a homomorphism
\[
\rho: (E\cup T\cup O)^+\to (S\cup E\cup Q\cup O)^*
\]
by $e_i\rho=e_i$ and $o_i\rho=o_i$ (for each $1\leq i\leq n$) and $t_{ij}\rho=\tau_{ij}$ (for each $\oijn$).  It follows that $\wb=\overline{w\rho}$ (i.e., $w\Phi=(w\rho)\Psi$) for all $w\in(E\cup T\cup O)^+$.  Note that $\im(\rho)$ is the subsemigroup of $(S\cup E\cup Q\cup O)^*$ generated by $E\cup O\cup\set{\tau_{ij}}{\oijn}$.  

Recall that $\sim$ is the congruence on $(E\cup T\cup O)^+$ generated by relations (R1--R17).

\ms
\begin{lemma}\label{lem:reduction1}
For any $u,v\in(E\cup T\cup O)^+$, $u\sim v \implies u\rho\approx v\rho$.
\end{lemma}

\pf First note that since $\sim$ is generated by relations (R1--R17), it suffices to prove the result for each 
relation $u=v$ from (R1--R17).
Since $w\rho\in(S\cup E\cup Q)^+$ for all $w\in(E\cup T)^+$, and since (as noted above) $\wb=\overline{w\rho}$ for all $w\in(E\cup T)^+$, it follows from Proposition \ref{prop:InPn} that the result is true for relations (R1--R10).  The result is also clearly true for relations (R11--R15) since these are precisely relations (R34), (R35) and (R39--41).  
It therefore remains to show that:
\bit
\item[(i)] $\tau_{ij}o_k\approx o_k\tau_{ij}$ for distinct $i,j$ and any $k$; and
\item[(ii)] $\tau_{ij}o_i\approx \tau_{ij}o_j \approx o_io_j$ for distinct $i,j$.
\eit
Beginning with (i), let $i,j,k\in\bn$ with $i\not=j$.  Put $l=k\Sb_{ij}^{-1}$.  Then $o_l\si_{ij}\approx \si_{ij}o_k$ and $\si_{ij}^{-1}o_l\approx o_k\si_{ij}^{-1}$ by Proposition \ref{prop:InPn} (and a simple diagrammatic check).  Together with (R42), we obtain
\[
\tau_{ij}o_k = \si_{ij}^{-1}q_i\si_{ij} o_k \approx \si_{ij}^{-1}q_io_l \si_{ij} \approx \si_{ij}^{-1}o_lq_i \si_{ij} \approx o_k \si_{ij}^{-1}q_i \si_{ij} =o_k \tau_{ij}.
\]
For (ii), let $i,j\in\bn$ with $i\not=j$.  Then, by Proposition \ref{prop:InPn}, (R43) and (R18), we have
\[
\tau_{ij}o_i = \si_{ij}^{-1}q_i\si_{ij} o_i \approx \si_{ij}^{-1}q_io_i\si_{ij}  \approx \si_{ij}^{-1}o_io_{i+1}\si_{ij} \approx o_i\si_{ij}^{-1}\si_{ij}o_j \approx o_io_j.
\]
A similar calculation yields $\tau_{ij}o_j\approx o_io_j$.  As noted above, this completes the proof. \epf

\nc{\ttb}{\overline{\tau}}

\ms
\begin{lemma}\label{lem:sxxs}
Let $1\leq i\leq n-1$ and $x\in E\cup O\cup\set{\tau_{ij}}{\oijn}$.  Then $s_ix$ and $xs_i$ are $\approx$-equivalent to an element of $\im(\rho)$.
\end{lemma}

\pf We consider separate cases, depending on whether $x$ belongs to $E$, $O$ or $\set{\tau_{ij}}{\oijn}$.

{\bf Case 1.}  Suppose first that $x\in E$.  One may check diagramatically that
\begin{align*}
\sb_i\eb_j = \eb_j\sb_i &= \eb_j\ttb_{i+1,j}\eb_{i+1}\ttb_{i,i+1}\eb_i\ttb_{ij}\eb_j &&\text{if $j\not\in\{i,i+1\}$}\\
\sb_i\eb_i = \eb_{i+1}\sb_i &= \eb_{i+1}\ttb_{i,i+1}\eb_i \\
\sb_i\eb_{i+1} = \eb_i\sb_i &= \eb_i\ttb_{i,i+1}\eb_{i+1} .
\intertext{It follows by Proposition \ref{prop:InPn} that}
s_ie_j \approx e_js_i &\approx e_j\tau_{i+1,j}e_{i+1}\tau_{i,i+1}e_i\tau_{ij}e_j &&\text{if $j\not\in\{i,i+1\}$}\\
s_ie_i \approx e_{i+1}s_i &\approx e_{i+1}\tau_{i,i+1}e_i \\
s_ie_{i+1} \approx e_is_i &\approx e_i\tau_{i,i+1}e_{i+1},
\end{align*}
completing the proof in this case.

{\bf Case 2.}  Next suppose $x=o_j\in O$, and put $k=j\sb_i$ (so also $j=k\sb_i$).  
Then 
 $o_js_i\approx s_io_k$
and 
$s_io_j\approx o_ks_i$, by Proposition \ref{prop:InPn}.  Together with (R40), we deduce that
\[
s_io_j \approx s_io_je_jo_j \approx o_k(s_ie_j)o_j
\AND
o_js_i \approx o_je_jo_js_i \approx o_j(e_js_i)o_k.
\]
By the previous case, $s_ie_j$ and $e_js_i$ are both $\approx$-equivalent to an element of $\im(\rho)$, so the proof is complete in this case also.

{\bf Case 3.}  Finally, suppose $x=\tau_{jk}$ for some $1\leq j<k\leq n$, and put $u=j\sb_i$ and $v=k\sb_i$.  By Proposition \ref{prop:InPn}, we have $\tau_{jk}s_i\approx s_i\tau_{uv}$ and $\tau_{jk}\approx \tau_{jk}e_j\tau_{jk}$.  It then follows that
\[
\tau_{jk}s_i \approx \tau_{jk}e_j\tau_{jk}s_i \approx \tau_{jk}(e_js_i)\tau_{uv},
\]
and, again, it follows from Case 1 that $\tau_{jk}s_i$ is $\approx$-equivalent to an element of $\im(\rho)$.  A similar calculation shows that this is the case also for $s_i\tau_{jk}$. \epf

\ms
\begin{lemma}\label{lem:reduction2}
Let $w\in(S\cup E\cup Q\cup O)^*\sm S^*$.  Then $w$ is $\approx$-equivalent to an element of $\im(\rho)$.
\end{lemma}

\pf Put $\Si=E\cup O\cup\set{\tau_{ij}}{\oijn}$, noting that $\im(\rho)=\la\Si\ra$ is the subsemigroup of $(S\cup E\cup Q\cup O)^*$ generated by $\Si$.  Since $Q\sub \Si$ (as $q_i=\tau_{i,i+1}$), it suffices to show that every element of $\la\Si\cup S\ra\sm S^*$ is $\approx$-equivalent to an element of $\im(\rho)$.  With this in mind, let $w\in\la\Si\cup S\ra\sm S^*$, and write $w=x_1\cdots x_k$, where $x_1,\ldots,x_k\in\Si\cup S$.  Denote by $l$ the number of factors $x_i$ that belong to~$S$.  We proceed by induction on $l$.  
If $l=0$, then we already have $w\in\la\Si\ra=\im(\rho)$, so suppose $l\geq1$.  Since $w\not\in S^*$, there exists $1\leq i\leq k-1$ such that either (i) $x_i\in S$ and $x_{i+1}\in\Si$, or (ii) $x_i\in\Si$ and $x_{i+1}\in S$.  In either case, Lemma \ref{lem:sxxs} tells us that $x_ix_{i+1}\approx u$ for some $u\in\im(\rho)=\la\Si\ra$.  But then $w\approx (x_1\cdots x_{i-1})u(x_{i+2}\cdots x_k)$, and we are done, after applying an induction hypothesis (noting that $(x_1\cdots x_{i-1})u(x_{i+2}\cdots x_k)$ has $l-1$ factors from $S$). \epf

With these preliminary results in place, we may now prove the first of the main results of this section.

\ms
\begin{thm}\label{thm:RPn}
The rook partition monoid, $\RP_n$, has monoid presentation 
\[
\pres{S\cup E\cup Q\cup O}{\text{\emph{(R18--R43)}}}
\]
via $\Psi$. 
\end{thm}

\pf It remains to show that $\ker\Psi\sub{\approx}$, so suppose $w_1,w_2\in(S\cup E\cup Q\cup O)^*$ are such that $\wb_1=\wb_2$.  If $\wb_1\in\S_n$, then $w_1,w_2\in S^*$, and $w_1\approx w_2$, using only relations (R18--R20).  For the remainder of the proof, suppose $\wb_1\not\in\S_n$.  It follows that $w_1,w_2\in(S\cup E\cup Q\cup O)^*\sm S^*$.  So, by Lemma \ref{lem:reduction2}, $w_1\approx u_1\rho$ and $w_2\approx u_2\rho$ for some $u_1,u_2\in(E\cup T\cup O)^+$.  We then have
\[
\ub_1=\overline{u_1\rho}=\wb_1=\wb_2=\overline{u_2\rho}=\ub_2,
\]
so that $u_1\sim u_2$, by Theorem \ref{thm:RPnSn}.  Lemma \ref{lem:reduction1} then gives $u_1\rho\approx u_2\rho$, so that $w_1\approx w_2$. \epf

\subsection{A presentation for $\RP_n$ on $n+2$ generators}\label{sect:RPn_n+2}

The presentation from Theorem \ref{thm:RPn} uses $4n-2$ generators.
In this section, we use \emph{Tietze transformations} to reduce the size of the generating set, thereby obtaining a presentation (Theorem \ref{thm:RPn2}) in terms of $S$ and three more generators. 
With this in mind, we rename $e=e_1$, $q=q_1$, $o=o_1$.  Define words
\begin{myalign}
\tag*{} c_i = s_1\cdots s_{i-1} \COMMA E_i=c_i^{-1}ec_i \COMMA O_i = c_i^{-1}oc_i &&&\text{for each $1\leq i\leq n$}\\
\tag*{} d_j = s_2\cdots s_j s_1\cdots s_{j-1} \COMMA Q_j = d_j^{-1}qd_j &&&\text{for each $1\leq j\leq n-1$.}
\end{myalign}
Note that $c_1=d_1=1$, $E_1=e$, $O_1=o$ and $Q_1=q$.  Diagramatically, one may check (in similar fashion to Figure~\ref{fig:tauij}) that 
\[
\eb_i = \Eb_i \COMMA \ob_i=\Ob_i \COMMA \qb_j=\Qb_j \qquad\text{for each $i,j$.}
\]
It follows from Theorem \ref{thm:RPn} that
\[
e_i\approx E_i \COMMA o_i\approx O_i \COMMA q_j\approx Q_j \qquad\text{for each $i,j$.}
\]
So we may transform the above presentation into $\pres{S\cup\{e,q,o\}}{\text{(R18--R43)$'$}}$, where each relation (R$k$)$'$ is obtained from (R$k$) by replacing each letter $e_ i,o_i,q_j$ by the words $E_i,O_i,Q_j$, respectively.  This presentation is via the map $\xi:(S\cup\{e,q,o\})^*\to\RP_n$ defined to be the restriction to $(S\cup\{e,q,o\})^*$ of $\Psi:(S\cup E\cup Q\cup O)^*\to\RP_n$.
Now consider the relations
\begin{align}
\tag{R44} s_i^2 &= 1 &&\text{for all $i$}\\
\tag{R45} s_is_j &= s_js_i &&\text{if $|i-j|>1$}\\
\tag{R46} s_is_js_i &= s_js_is_j &&\text{if $|i-j|=1$}\\
\tag{R47} e^2 = e &= e qe = eoe\\
\tag{R48} q^2=q=qe q&=qs_1=s_1q\\
\tag{R49} e s_i &= s_i e &&\text{if $i\geq2$}\\
\tag{R50} q s_i &= s_i q &&\text{if $i\geq3$}\\
\tag{R51} s_1e s_1e = e &s_1e s_1 = e s_1e\\
\tag{R52} q s_2 q s_2 &= s_2 q s_2 q \\
\tag{R53} q (s_2s_1s_3s_2) q (s_2s_1s_3s_2) &= (s_2s_1s_3s_2) q (s_2s_1s_3s_2) q\\
\tag{R54} q(s_2s_1e s_1s_2) &= (s_2s_1e s_1s_2)q \\
\tag{R55} o^2=o&=oeo \\
\tag{R56} os_i&=s_io &&\text{if $i\geq2$} \\
\tag{R57} os_1os_1=s_1os_1o&=os_1o=oq=qo \\
\tag{R58} es_1os_1&=s_1os_1e \\
\tag{R59} q(s_2s_1o s_1s_2) &= (s_2s_1o s_1s_2)q.
\end{align}
It is easy to check, diagramatically, that each of relations (R44--R59) hold in $\RP_n$ when the words are replaced by their images under $\xi$.  So we may add them to the presentation, obtaining the presentation $\pres{S\cup\{e,q,o\}}{\text{(R18--R43)$'$, (R44--R59)}}$.  Our goal is to show that the relations (R18--R43)$'$ may be removed; see Theorem \ref{thm:RPn2}.
With this in mind, write $\ssim$ for the congruence on $(S\cup\{e,q,o\})^*$ generated by (R44--R59).  For $w\in (S\cup\{e,q,o\})^*$, write $\wb=w\xi$.

\ms
\begin{prop}\label{prop:InPn2}
If $u,v\in(S\cup \{e,q\})^*$ or $u,v\in(S\cup \{o\})^*$, then $\ub=\vb\implies u\ssim v$.
\end{prop}

\pf The proof is similar to that of Proposition \ref{prop:InPn}, since (R44-R54) include defining relations for 
$\P_n=\la S\xi\cup\{\eb,\qb\}\ra$ 
(see \cite[Theorem 32]{JEgrpm}), while (R44--R45) and (R55--R57) include defining relations for 
$\R_n=\la S\xi\cup\{\ob\}\ra$ 
(see \cite{Popova} or \cite{Fernandes2001}). \epf

In particular, we may remove any of the relations (R$k$)$'$ involving only words over $S\cup\{e,q\}$ or only words over $S\cup\{o\}$.  In this way, we may remove relations (R18--R38)$'$.
The next two results follow immediately from Proposition \ref{prop:InPn2} (and simple diagrammatic checks).

\ms
\begin{cor}\label{cor:Eiw}
Let $w\in S^*$ and $i\in\bn$.  Then $w^{-1}E_iw\ssim E_{i\wb}$ and $w^{-1}O_iw\ssim O_{i\wb}$. \epfres
\end{cor}

\ms
\begin{cor}\label{cor:Qiw}
Let $w\in S^*$ and $1\leq i\leq n-1$, with $(i+1)\wb=i\wb+1$.  Then ${w^{-1}Q_iw\ssim Q_{i\wb}}$.~\epfres
\end{cor}

\ms
\begin{thm}\label{thm:RPn2}
The rook partition monoid, $\RP_n$, has monoid presentation 
\[
\pres{S\cup\{e,q,o\}}{\text{\emph{(R44--R59)}}}
\]
via $\xi$. 
\end{thm}

\pf So far, we have transformed the presentation into $\pres{S\cup\{e,q,o\}}{\text{(R39--R43)$'$, (R44--R59)}}$.  It remains to show that relations (R39--R43)$'$ may be removed.

{\bf (R39)$'$.}  We need to show that $O_iE_j\ssim E_jO_i$ if $i\not=j$.  Let $w\in S^*$ be such that $1\wb=j$ and $2\wb=i$.  Then
\[
O_iE_j \ssim w^{-1}O_2ww^{-1}E_1w \ssim w^{-1}O_2E_1w \ssim w^{-1}E_1O_2w \ssim w^{-1}E_1ww^{-1}O_2w \ssim E_jO_i,
\]
by Corollary \ref{cor:Eiw}, (R44) and (R58), the latter of which says ``$E_1O_2=O_2E_1$''.

{\bf (R40)$'$ and (R41)$'$.}  Here we have $O_iE_iO_i = c_i^{-1}oc_ic_i^{-1}ec_ic_i^{-1}oc_i \ssim c_i^{-1}oeoc_i \ssim c_i^{-1}oc_i = O_i$, by (R44) and (R55).  An almost identical calculation (using (R47) instead of (R55)) deals with relation (R41)$'$.

{\bf (R42)$'$.}  We must show that $Q_iO_j\ssim O_jQ_i$ for any $i,j$.  Suppose first that $j\not\in\{i,i+1\}$.  Choose $w\in S^*$ such that $1\wb=i$, $2\wb=i+1$ and $3\wb=j$.  Then
\[
Q_iO_j \ssim w^{-1}Q_1w w^{-1}O_3w \ssim w^{-1}Q_1O_3w \ssim w^{-1}O_3Q_1w \ssim w^{-1}O_3ww^{-1}Q_1w \ssim O_jQ_i,
\]
by Corollaries \ref{cor:Eiw} and \ref{cor:Qiw}, and relations (R44) and (R59), the latter of which says ``$Q_1O_3=O_3Q_1$''.  Next, note that for any $u\in S^*$ with $1\ub=i$ and $2\ub=i+1$,
$Q_iO_i \ssim u^{-1} qo u \ssim u^{-1}oqu \ssim O_iQ_i$, by Corollaries \ref{cor:Eiw} and \ref{cor:Qiw}, and relations (R44) and (R57).  A similar calculation gives $Q_iO_{i+1}\ssim O_{i+1}Q_i$.

{\bf (R43)$'$.}  Again, let $u\in S^*$ be such that $1\ub=i$ and $2\ub=i+1$.  Then $Q_iO_i \ssim u^{-1}Q_1O_1u \ssim u^{-1}O_1O_2u \ssim O_iO_{i+1}$, by Corollaries \ref{cor:Eiw} and \ref{cor:Qiw}, and relations (R44) and (R57), the relevant part of the latter of which says ``$O_1Q_1=O_1O_2$''.  Finally, using Proposition \ref{prop:InPn2} and the previous calculation, we also have $Q_iO_{i+1} = Q_is_iO_is_i  \ssim Q_iO_is_i \ssim O_i O_{i+1}s_i \ssim O_i O_{i+1}$.
This completes the proof. \epf

\subsection{A presentation for $\RP_n$ on $5$ generators}\label{sect:RPn_5}

We continue to reduce the presentation $\pres{S\cup\{e,q,o\}}{\text{(R44--R59)}}$ from Theorem \ref{thm:RPn2} in order to obtain a presentation with the minimal number of generators, making use of a $2$-generator presentation for $\S_n$ \cite{Moo}.  In fact, the generating set $S\cup\{e,q,o\}$ is already of minimal size when $n=2$ (see Theorem~\ref{thm:rankRPn} and Remark \ref{rem:rankRPn}), so we will assume $n\geq3$ for the rest of this section.  

We now rename $s=s_1$, and add the new generator $c$, along with the relation $c=s_1\cdots s_{n-1}$.  It is easy to check, diagrammatically, that $\sb_{i+1}=\cb^i\;\! \sb\ \cb^{-i}=\cb^i\;\!\sb\ \cb^{n-i}$ for all $2\leq i\leq n-2$; here we write $\cb=\sb_1\cdots\sb_{n-1}$.  It follows that we may remove the generators $s_2,\ldots,s_{n-1}$ from the presentation, and replace their every occurrence in the relations by the words $S_{i+1}=c^isc^{n-i}$ (for each $1\leq i\leq n-2$).  (We also define $S_1=s$.)  The result of doing this to relation (R$k$) will be denoted (R$k$)$'$.  So we have transformed the presentation to
\[
\pres{s,c,e,q,o}{\text{(R44--R59)$'$},\ c=S_1\cdots S_{n-1}}.
\]
This presentation is via the map
\[
\zeta:\{s,c,e,q,o\}^*\to\RP_n:x\mt\xb.
\]
Now consider the relations
\begin{gather}
\tag{R60} c^n = (sc)^{n-1} = s^2 = (c^isc^{n-i}s)^2 =1 \qquad\text{for all $2\leq   i\leq \tfrac{n}{2}$}\\
\tag{R61} e^2=e=e qe = eoe = sce c^{n-1}s = csc^{n-1} e csc^{n-1}\\
\tag{R62} q^2=q=qeq = qs=sq= c^2sc^{n-2} q c^2sc^{n-2} =   c^{n-1}scsc^{n-1}qcsc^{n-1}sc \\
\tag{R63} se se = e s e s = e s e           \\
\tag{R64}  qcqc^{n-1} = cqc^{n-1}q          \\
\tag{R65}  qc^2qc^{n-2} = c^2qc^{n-2}q          \\
\tag{R66} qc^2e c^{n-2} = c^2e c^{n-2}q          \\
\tag{R67} o^2=o=oeo = sco c^{n-1}s = csc^{n-1} o csc^{n-1}\\
\tag{R68} so so = o s o s = o s o  =oq=qo         \\
\tag{R69} qc^2o c^{n-2} = c^2o c^{n-2}q          \\
\tag{R70} esos = sose          .
\end{gather}

\ms
\begin{thm}\label{thm:RPn3}
The rook partition monoid, $\RP_n$, has monoid presentation 
\[
\pres{s,c,e,q,o}{\text{\emph{(R60--R70)}}}
\]
via $\zeta$. 
\end{thm}

\pf Write $\diamond$ for the congruence on $\{s,c,e,q,o\}^*$ generated by relations (R60--R70).  We already know that $\zeta$ is surjective, and it is easy to check that $\diamond\sub\ker\zeta$, so we may add relations (R60--R70) to obtain the presentation
\[
\pres{s,c,e,q,o}{\text{(R60--R70)},\ \text{(R44--R59)$'$},\ c=S_1\cdots S_{n-1}}.
\]
It remains to show that we can remove all relations apart from (R60--R70).
Relations (R60--R66) contain defining relations for $\P_n=\la\sb,\cb,\eb,\qb\ra$ (see \cite[Theorem 41]{JEgrpm}), so we may remove $c=S_1\cdots S_{n-1}$ and all relations from (R44--R59)$'$ with no occurrence of the letter $o$.  This leaves (R55--R59)$'$ and the relation $e=eoe$; the latter is part of (R61), so may be removed.
Next note that (R60), (R67) and (R68) contain defining relations for $\R_n=\la\sb,\cb,\ob\ra$ (see \cite{Popova} or \cite{Fernandes2001}), so we may remove all relations from (R55--R59)$'$ with no occurrence of the letters $e,q$.  This leaves (R58)$'$, (R59)$'$, and the relations
\[
o=oeo \AND oq=qo=oso.
\]
These last relations are already part of (R67) and (R68), so we are only left with (R58)$'$ and (R59)$'$.  The former is just (R70).  For (R59)$'$, first observe that $\sb_2\sb_1\ob\ \sb_1\sb_2=\cb^2\;\!\ob\ \cb^{n-2}$, so that $S_2S_1oS_1S_2\mathrel{\diamond} c^2oc^{n-2}$, by the above-mentioned fact that (R60--R70) contains defining relations for $\R_n=\la\sb,\cb,\ob\ra$.  Together with (R69), it then follows that
$qS_2S_1oS_1S_2 \mathrel{\diamond} qc^2oc^{n-2}  \mathrel{\diamond} c^2oc^{n-2} q \mathrel{\diamond} S_2S_1oS_1S_2q.$ \epf

\ms
\begin{rem}
The presentation in Theorem \ref{thm:RPn3} includes defining relations for the symmetric group~$\S_n$ (i.e., relations (R60)), plus a fixed (i.e., independent of $n$) number of extra relations: 26 such extra relations, to be precise.  We make no claim that this is the minimal number of additional relations required.
\end{rem}

We conclude this section by showing that the generating set in the previous result has the minimal possible size.

\ms
\begin{thm}\label{thm:rankRPn}
For $n\geq3$, $\rank(\RP_n)=5$.
\end{thm}

\pf By Theorem \ref{thm:RPn3}, it suffices to show that any generating set for $\RP_n$ has size at least $5$.  So suppose $\RP_n=\la\Si\ra$.  Since $\RPnSn$ is an ideal of $\RP_n$, it follows that $\Si\cap\S_n$ is a generating set for~$\S_n$.  Since $\rank(\S_n)=2$, it follows that $|\Si\cap\S_n|\geq2$.  It therefore remains to show that $|\Si\sm\S_n|\geq3$.

Consider an expression $\eb_1=\al_1\cdots\al_k$, where $\al_1,\ldots,\al_k\in\Si$.  Let $1\leq l\leq k$ be minimal so that $\al_l\not\in\S_n$.
Let $i=1(\al_1\cdots\al_{l-1})$.  Then 
\[
\eb_i = (\al_1\cdots\al_{l-1})^{-1}\eb_1(\al_1\cdots\al_{l-1}) = \al_l(\al_{l+1}\cdots\al_k)(\al_1\cdots\al_{l-1}).
\]
As in the proof of Theorem \ref{thm:rankRPnSn}, it follows that $\al_l$ has domain $\{i\}^c$ and the block $\{i\}$.  Similarly, it can be shown that $\Si\sm\S_n$ contains: an element with domain $\{j\}^c$ and the rook dot $j$, for some $j\in\bn$; and an element with domain $\bn$ and non-trivial kernel. \epf

\ms
\begin{rem}\label{rem:rankRPn}
The second paragraph of the previous proof works for $n=2$ as well, showing that $\rank(\RP_2)=\rank(\S_2)+3=4$.  It is easy to check that $\rank(\RP_n)=1,3$ for $n=0,1$.
\end{rem}

\section{Rook partition algebras}\label{sect:algebras}

Recall from \cite{Wilcox2007} that the \emph{partition algebras} are \emph{twisted semigroup algebras} of the partition monoids.  As noted in the introduction, this understanding has led to a rich flow of information between the theories of diagram \emph{semigroups} and diagram \emph{algebras}: see for example \cite{Wilcox2007,JEgrpm,JEpnsn,DEG2015,EastGray,DEEFHHL1,DEEFHHL2}.  In this section, we show how the results of previous sections lead to (algebra) presentations for the \emph{rook partition algebras} and their singular ideals.

Let $F$ be a commutative ring with $1$, and let $S$ be a semigroup.  Recall that a \emph{twisting} from $S$ to~$F$ is a map $\tau:S\times S\to F$ satisfying $\tau(a,b)\tau(ab,c)=\tau(a,bc)\tau(b,c)$ for all $a,b,c\in S$.
Given such a twisting, the \emph{twisted semigroup algebra} of $S$ over $F$ with respect to $\tau$, denoted $F^\tau[S]$, is defined to be the set of all finite formal $F$-linear combinations over $S$, with associative operation $\star$ defined on basis elements (and then extended linearly) by $a\star b=\tau(a,b)ab$ for each $a,b\in S$ (where $ab$ is the product of $a$ and $b$ in $S$).  When the twisting is trivial (i.e., $\tau(a,b)=1\in F$ for all $a,b\in S$), $F^\tau[S]=F[S]$ is just the (ordinary) semigroup algebra.

Given a fixed element $\de\in F$, one may define a twisting from $\RP_n$ to $F$ as follows.  For $\al,\be\in\RP_n$, let $m(\al,\be)$ denote the number of connected components in the product graph $\Ga(\al,\be)$ that involve only black vertices in the middle layer (i.e., double-dashed non-rook vertices).  One may show that
\[
m(\al,\be)+m(\al\be,\ga)=m(\al,\be\ga)+m(\be,\ga) \qquad\text{for all $\al,\be,\ga\in\RP_n$.}
\]
As a result, one may then define a twisting
\[
\tau:\RP_n\times\RP_n\to F \qquad\text{by}\qquad \tau(\al,\be)=\de^{m(\al,\be)} \qquad\text{for all $\al,\be\in\RP_n$,}
\]
and the resulting twisted semigroup algebra $F^\tau[\RP_n]$ is called the \emph{rook partition algebra} \cite{Grood06}.  (The reason that white vertices (i.e., rook dots) do not figure in the count of $m(\al,\be)$ is due to the exact nature of the above-mentioned embedding of $\RP_n$ in $\P_{n+1}$; see \cite{Grood06} for more details.)  Note also that~$\tau$ restricts to a twisting from $\RPnSn$ to $F$, so we may consider the \emph{singular rook partition algebra}, $F^\tau[\RPnSn]$.

A general result from \cite{JEgrpm} shows how a (monoid or semigroup) presentation for $S$ leads to an (algebra) presentation for $F^\tau[S]$ in the case that the image of the twisting $\tau$ lies in $G(F)$, the group of units of $F$.  The reader is referred to \cite[Section 6]{JEgrpm} for full details.  The next result follows immediately from \cite[Theorem 44]{JEgrpm} and earlier results of the current paper (as stated below).

\ms
\begin{thm}\label{thm:algebras}
Suppose $F$ is a commutative ring with identity, and let $\de\in G(F)$ be a unit in $F$.
\bit
\itemit{i} An algebra presentation for the rook partition algebra, $F^\tau[\RP_n]$, may be obtained from any of the above monoid presentations for $\RP_n$ by:
\begin{itemize}
\itemit{a} changing the relations $e_i^2=e_i$ to $e_i^2=\de e_i$ in Theorem \ref{thm:RPn}; or
\itemit{b} changing the relations $e^2=e$ to $e^2=\de e$ in Theorem \ref{thm:RPn2} or Theorem \ref{thm:RPn3}.
\end{itemize}
\itemit{ii} An algebra presentation for the singular rook partition algebra, $F^\tau[\RPnSn]$, may be obtained from the semigroup presentation for $\RPnSn$ in Theorem \ref{thm:RPnSn} by changing the relations $e_i^2=e_i$ to $e_i^2=\de e_i$. \epfres
\eit
\end{thm}

\footnotesize
\def\bibspacing{-1.1pt}
\bibliography{biblio}

\begin{thebibliography}{10}

\bibitem{ACHLV2015}
K.~Auinger, Yuzhu Chen, Xun Hu, Yanfeng Luo, and M.~V. Volkov.
\newblock The finite basis problem for {K}auffman monoids.
\newblock {\em Algebra Universalis}, 74(3-4):333--350, 2015.

\bibitem{Auinger2014}
Karl Auinger.
\newblock Pseudovarieties generated by {B}rauer type monoids.
\newblock {\em Forum Math.}, 26(1):1--24, 2014.

\bibitem{ADV2012_2}
Karl Auinger, Igor Dolinka, and Mikhail~V. Volkov.
\newblock Equational theories of semigroups with involution.
\newblock {\em J. Algebra}, 369:203--225, 2012.

\bibitem{ADV2012}
Karl Auinger, Igor Dolinka, and Mikhail~V. Volkov.
\newblock Matrix identities involving multiplication and transposition.
\newblock {\em J. Eur. Math. Soc. (JEMS)}, 14(3):937--969, 2012.

\bibitem{BH2014}
Georgia Benkart and Tom Halverson.
\newblock Motzkin algebras.
\newblock {\em European J. Combin.}, 36:473--502, 2014.

\bibitem{Bloss2003}
Matthew Bloss.
\newblock {$G$}-colored partition algebras as centralizer algebras of wreath
  products.
\newblock {\em J. Algebra}, 265(2):690--710, 2003.

\bibitem{BDP2002}
Mirjana Borisavljevi{\'c}, Kosta Do{\v{s}}en, and Zoran Petri{\'c}.
\newblock Kauffman monoids.
\newblock {\em J. Knot Theory Ramifications}, 11(2):127--143, 2002.

\bibitem{Brauer1937}
Richard Brauer.
\newblock On algebras which are connected with the semisimple continuous
  groups.
\newblock {\em Ann. of Math. (2)}, 38(4):857--872, 1937.

\bibitem{Cameron1999}
Peter~J. Cameron.
\newblock {\em Permutation groups}, volume~45 of {\em London Mathematical
  Society Student Texts}.
\newblock Cambridge University Press, Cambridge, 1999.

\bibitem{CGM2003}
Anton Cox, John Graham, and Paul Martin.
\newblock The blob algebra in positive characteristic.
\newblock {\em J. Algebra}, 266(2):584--635, 2003.

\bibitem{DO2014}
Zajj Daugherty and Rosa Orellana.
\newblock The quasi-partition algebra.
\newblock {\em J. Algebra}, 404:124--151, 2014.

\bibitem{DHP2003}
Momar Dieng, Tom Halverson, and Vahe Poladian.
\newblock Character formulas for {$q$}-rook monoid algebras.
\newblock {\em J. Algebraic Combin.}, 17(2):99--123, 2003.

\bibitem{DE2015}
Igor Dolinka and James East.
\newblock Twisted {B}rauer monoids.
\newblock {\em Preprint}, 2015, {\tt arXiv:1510.08666}.

\bibitem{DE2016}
Igor Dolinka and James East.
\newblock The idempotent generated subsemigroup of the {K}auffman monoid.
\newblock {\em Preprint}, 2016, {\tt arXiv:1602.01157}.

\bibitem{DEEFHHL1}
Igor Dolinka, James East, Athanasios Evangelou, Des FitzGerald, Nicholas Ham,
  James Hyde, and Nicholas Loughlin.
\newblock Enumeration of idempotents in diagram semigroups and algebras.
\newblock {\em J. Combin. Theory Ser. A}, 131:119--152, 2015.

\bibitem{DEEFHHL2}
Igor Dolinka, James East, Athanasios Evangelou, Des FitzGerald, Nicholas Ham,
  James Hyde, and Nicholas Loughlin.
\newblock Idempotent statistics of the {M}otzkin and {J}ones monoids.
\newblock {\em Preprint}, 2015, {\tt arXiv:1507.04838}.

\bibitem{DEG2015}
Igor Dolinka, James East, and Robert~D. Gray.
\newblock Motzkin monoids and partial {B}rauer monoids.
\newblock {\em Preprint}, 2015, {\tt arXiv:1512.02279}.

\bibitem{JEcais}
James East.
\newblock Cellular algebras and inverse semigroups.
\newblock {\em J. Algebra}, 296(2):505--519, 2006.

\bibitem{JEboppp}
James East.
\newblock Braids and order-preserving partial permutations.
\newblock {\em J. Knot Theory Ramifications}, 19(8):1025--1049, 2010.

\bibitem{JEgrpm}
James East.
\newblock Generators and relations for partition monoids and algebras.
\newblock {\em J. Algebra}, 339:1--26, 2011.

\bibitem{JEpnsn}
James East.
\newblock On the singular part of the partition monoid.
\newblock {\em Internat. J. Algebra Comput.}, 21(1-2):147--178, 2011.

\bibitem{JEinsn2}
James East.
\newblock A symmetrical presentation for the singular part of the symmetric
  inverse monoid.
\newblock {\em Algebra Universalis}, 74(3-4):207--228, 2015.

\bibitem{EF}
James East and D.~G. FitzGerald.
\newblock The semigroup generated by the idempotents of a partition monoid.
\newblock {\em J. Algebra}, 372:108--133, 2012.

\bibitem{EastGray}
James East and R.~D. Gray.
\newblock Diagram monoids and {G}raham--{H}oughton graphs: idempotents and
  generating sets of ideals.
\newblock {\em Preprint}, 2014, {\tt arXiv:1404.2359}.

\bibitem{Enyang2013_1}
John Enyang.
\newblock Jucys-{M}urphy elements and a presentation for partition algebras.
\newblock {\em J. Algebraic Combin.}, 37(3):401--454, 2013.

\bibitem{Enyang2013_2}
John Enyang.
\newblock A seminormal form for partition algebras.
\newblock {\em J. Combin. Theory Ser. A}, 120(7):1737--1785, 2013.

\bibitem{Fernandes2001}
V{\'{\i}}tor~H. Fernandes.
\newblock Presentations for some monoids of partial transformations on a finite
  chain: a survey.
\newblock In {\em Semigroups, algorithms, automata and languages ({C}oimbra,
  2001)}, pages 363--378. World Sci. Publ., River Edge, NJ, 2002.

\bibitem{FL2011}
D.~G. FitzGerald and Kwok~Wai Lau.
\newblock On the partition monoid and some related semigroups.
\newblock {\em Bull. Aust. Math. Soc.}, 83(2):273--288, 2011.

\bibitem{FL1998}
D.~G. FitzGerald and Jonathan Leech.
\newblock Dual symmetric inverse monoids and representation theory.
\newblock {\em J. Austral. Math. Soc. Ser. A}, 64(3):345--367, 1998.

\bibitem{Gilbert2006}
N.~D. Gilbert.
\newblock Presentations of the inverse braid monoid.
\newblock {\em J. Knot Theory Ramifications}, 15(5):571--588, 2006.

\bibitem{GL1996}
J.~J. Graham and G.~I. Lehrer.
\newblock Cellular algebras.
\newblock {\em Invent. Math.}, 123(1):1--34, 1996.

\bibitem{Grood2002}
Cheryl Grood.
\newblock A {S}pecht module analog for the rook monoid.
\newblock {\em Electron. J. Combin.}, 9(1):Research Paper 2, 10 pp.
  (electronic), 2002.

\bibitem{Grood06}
Cheryl Grood.
\newblock The rook partition algebra.
\newblock {\em J. Combin. Theory Ser. A}, 113(2):325--351, 2006.

\bibitem{Halverson2004}
Tom Halverson.
\newblock Representations of the {$q$}-rook monoid.
\newblock {\em J. Algebra}, 273(1):227--251, 2004.

\bibitem{HD2014}
Tom Halverson and Elise delMas.
\newblock Representations of the {R}ook-{B}rauer algebra.
\newblock {\em Comm. Algebra}, 42(1):423--443, 2014.

\bibitem{HR2001}
Tom Halverson and Arun Ram.
\newblock {$q$}-rook monoid algebras, {H}ecke algebras, and {S}chur-{W}eyl
  duality.
\newblock {\em Zap. Nauchn. Sem. S.-Peterburg. Otdel. Mat. Inst. Steklov.
  (POMI)}, 283(Teor. Predst. Din. Sist. Komb. i Algoritm. Metody. 6):224--250,
  262--263, 2001.

\bibitem{HR2005}
Tom Halverson and Arun Ram.
\newblock Partition algebras.
\newblock {\em European J. Combin.}, 26(6):869--921, 2005.

\bibitem{Hig}
Peter~M. Higgins.
\newblock {\em Techniques of semigroup theory}.
\newblock Oxford Science Publications. The Clarendon Press, Oxford University
  Press, New York, 1992.

\bibitem{Howie1966}
J.~M. Howie.
\newblock The subsemigroup generated by the idempotents of a full
  transformation semigroup.
\newblock {\em J. London Math. Soc.}, 41:707--716, 1966.

\bibitem{Howie1978}
J.~M. Howie.
\newblock Idempotent generators in finite full transformation semigroups.
\newblock {\em Proc. Roy. Soc. Edinburgh Sect. A}, 81(3-4):317--323, 1978.

\bibitem{Howie}
John~M. Howie.
\newblock {\em Fundamentals of semigroup theory}, volume~12 of {\em London
  Mathematical Society Monographs. New Series}.
\newblock The Clarendon Press, Oxford University Press, New York, 1995.
\newblock Oxford Science Publications.

\bibitem{Jones1987}
V.~F.~R. Jones.
\newblock Hecke algebra representations of braid groups and link polynomials.
\newblock {\em Ann. of Math. (2)}, 126(2):335--388, 1987.

\bibitem{Jones1994_2}
V.~F.~R. Jones.
\newblock The {P}otts model and the symmetric group.
\newblock In {\em Subfactors ({K}yuzeso, 1993)}, pages 259--267. World Sci.
  Publ., River Edge, NJ, 1994.

\bibitem{Kauffman1990}
Louis~H. Kauffman.
\newblock An invariant of regular isotopy.
\newblock {\em Trans. Amer. Math. Soc.}, 318(2):417--471, 1990.

\bibitem{Kennedy2007}
A.~Joseph Kennedy.
\newblock Class partition algebras as centralizer algebras of wreath products.
\newblock {\em Comm. Algebra}, 35(1):145--170, 2007.

\bibitem{KM2013}
A.~Joseph Kennedy and G.~Muniasamy.
\newblock Rook version of colored partition algebras.
\newblock {\em Bull. Math. Sci.}, 3(1):1--17, 2013.

\bibitem{KM2011}
Ganna Kudryavtseva and Victor Maltcev.
\newblock Two generalisations of the symmetric inverse semigroups.
\newblock {\em Publ. Math. Debrecen}, 78(2):253--282, 2011.

\bibitem{KMU2015}
Ganna Kudryavtseva, Victor Maltcev, and Abdullahi Umar.
\newblock Presentation for the partial dual symmetric inverse monoid.
\newblock {\em Comm. Algebra}, 43(4):1621--1639, 2015.

\bibitem{KM2006}
Ganna Kudryavtseva and Volodymyr Mazorchuk.
\newblock On presentations of {B}rauer-type monoids.
\newblock {\em Cent. Eur. J. Math.}, 4(3):413--434 (electronic), 2006.

\bibitem{KM2008}
Ganna Kudryavtseva and Volodymyr Mazorchuk.
\newblock Schur-{W}eyl dualities for symmetric inverse semigroups.
\newblock {\em J. Pure Appl. Algebra}, 212(8):1987--1995, 2008.

\bibitem{LF2006}
Kwok~Wai Lau and D.~G. FitzGerald.
\newblock Ideal structure of the {K}auffman and related monoids.
\newblock {\em Comm. Algebra}, 34(7):2617--2629, 2006.

\bibitem{Lawson1998}
Mark~V. Lawson.
\newblock {\em Inverse semigroups}.
\newblock World Scientific Publishing Co., Inc., River Edge, NJ, 1998.
\newblock The theory of partial symmetries.

\bibitem{LZ2015}
Gustav Lehrer and R.~B. Zhang.
\newblock The {B}rauer category and invariant theory.
\newblock {\em J. Eur. Math. Soc. (JEMS)}, 17(9):2311--2351, 2015.

\bibitem{LZ2012}
Gustav Lehrer and Ruibin Zhang.
\newblock The second fundamental theorem of invariant theory for the orthogonal
  group.
\newblock {\em Ann. of Math. (2)}, 176(3):2031--2054, 2012.

\bibitem{Lipscombe1996}
Stephen Lipscomb.
\newblock {\em Symmetric inverse semigroups}, volume~46 of {\em Mathematical
  Surveys and Monographs}.
\newblock American Mathematical Society, Providence, RI, 1996.

\bibitem{Maltcev2007}
Victor Maltcev and Volodymyr Mazorchuk.
\newblock Presentation of the singular part of the {B}rauer monoid.
\newblock {\em Math. Bohem.}, 132(3):297--323, 2007.

\bibitem{Martin2000}
P.~P. Martin.
\newblock The partition algebra and the {P}otts model transfer matrix spectrum
  in high dimensions.
\newblock {\em J. Phys. A}, 33(19):3669--3695, 2000.

\bibitem{Martin1994}
Paul Martin.
\newblock Temperley-{L}ieb algebras for nonplanar statistical mechanics---the
  partition algebra construction.
\newblock {\em J. Knot Theory Ramifications}, 3(1):51--82, 1994.

\bibitem{Martin1996}
Paul Martin.
\newblock The structure of the partition algebras.
\newblock {\em J. Algebra}, 183(2):319--358, 1996.

\bibitem{MarMaz2014}
Paul Martin and Volodymyr Mazorchuk.
\newblock On the representation theory of partial {B}rauer algebras.
\newblock {\em Q. J. Math.}, 65(1):225--247, 2014.

\bibitem{MS1994}
Paul Martin and Hubert Saleur.
\newblock The blob algebra and the periodic {T}emperley-{L}ieb algebra.
\newblock {\em Lett. Math. Phys.}, 30(3):189--206, 1994.

\bibitem{MW2000}
Paul~P. Martin and David Woodcock.
\newblock On the structure of the blob algebra.
\newblock {\em J. Algebra}, 225(2):957--988, 2000.

\bibitem{Maz1998}
Volodymyr Mazorchuk.
\newblock On the structure of {B}rauer semigroup and its partial analogue.
\newblock {\em Problems in Algebra}, 13:29--45, 1998.

\bibitem{Maz2002}
Volodymyr Mazorchuk.
\newblock Endomorphisms of {$\mathfrak B_n$, $\mathcal P \mathfrak B_n$}, and
  {$\mathfrak C_n$}.
\newblock {\em Comm. Algebra}, 30(7):3489--3513, 2002.

\bibitem{Moo}
Eliakim~Hastings Moore.
\newblock Concerning the abstract groups of order $k!$ and $\tfrac{1}{2}k!$
  holohedrically isomorphic with the symmetric and the alternating
  substitution-groups on $k$ letters.
\newblock {\em Proc. London Math. Soc.}, 28(1):357--366, 1897.

\bibitem{NS1978}
T.~E. Nordahl and H.~E. Scheiblich.
\newblock Regular {$\ast $}-semigroups.
\newblock {\em Semigroup Forum}, 16(3):369--377, 1978.

\bibitem{Orellana2007}
Rosa~C. Orellana.
\newblock On partition algebras for complex reflection groups.
\newblock {\em J. Algebra}, 313(2):590--616, 2007.

\bibitem{Paget2006}
Rowena Paget.
\newblock Representation theory of {$q$}-rook monoid algebras.
\newblock {\em J. Algebraic Combin.}, 24(3):239--252, 2006.

\bibitem{PK2004}
M.~Parvathi and A.~Joseph Kennedy.
\newblock {$G$}-vertex colored partition algebras as centralizer algebras of
  direct products.
\newblock {\em Comm. Algebra}, 32(11):4337--4361, 2004.

\bibitem{Popova}
L.~M. Popova.
\newblock Defining relations in some semigroups of partial transformations of a
  finite set (in {R}ussian).
\newblock {\em Uchenye Zap. Leningrad Gos. Ped. Inst.}, 218:191--212, 1961.

\bibitem{PHY2013}
Eliezer Posner, Kris Hatch, and Megan Ly.
\newblock Presentation of the {M}otzkin monoid.
\newblock {\em Preprint}, 2013, {\tt arXiv:1301.4518}.

\bibitem{Solomon2002}
Louis Solomon.
\newblock Representations of the rook monoid.
\newblock {\em J. Algebra}, 256(2):309--342, 2002.

\bibitem{TL1971}
H.~N.~V. Temperley and E.~H. Lieb.
\newblock Relations between the ``percolation'' and ``colouring'' problem and
  other graph-theoretical problems associated with regular planar lattices:
  some exact results for the ``percolation'' problem.
\newblock {\em Proc. Roy. Soc. London Ser. A}, 322(1549):251--280, 1971.

\bibitem{Weyl_book}
Hermann Weyl.
\newblock {\em The classical groups}.
\newblock Princeton Landmarks in Mathematics. Princeton University Press,
  Princeton, NJ, 1997.
\newblock Their invariants and representations, Fifteenth printing, Princeton
  Paperbacks.

\bibitem{Wilcox2007}
Stewart Wilcox.
\newblock Cellularity of diagram algebras as twisted semigroup algebras.
\newblock {\em J. Algebra}, 309(1):10--31, 2007.

\end{thebibliography}
\bibliographystyle{plain}
\end{document}